\documentclass[11pt,reqno]{amsart}
\usepackage{amssymb}
\usepackage{amsmath}
\usepackage{amsfonts}
\usepackage{amsthm,latexsym,color}

\usepackage{pifont}
\usepackage{verbatim}
\usepackage{amsfonts}
\usepackage{hyperref}

\font\TenEns=msbm10 \font\SevenEns=msbm7 \font\FiveEns=msbm5
\newfam\Ensfam

\textfont\Ensfam=\TenEns \scriptfont\Ensfam=\SevenEns
\scriptscriptfont\Ensfam=\FiveEns

\def\R {\mathbb R}
\def\C {\mathbb C}

\DeclareMathOperator{\e}{e}

\newtheorem{theorem}{Theorem}
\theoremstyle{plain}

\newtheorem{definition}{Definition}

\newtheorem{lemma}{Lemma}[section]

\newtheorem{proposition}{Proposition}[section]
\newtheorem{remark}{Remark}[section]

\numberwithin{equation}{section}
\numberwithin{equation}{section}
\title[Damped Inhomogeneous Nonlinear Schr\"odinger Equation]{Asymptotic behavior  and life-span estimates for the damped Inhomogeneous Nonlinear Schr\"odinger Equation}
\author[L. Aloui]{}
\email{lassaad.aloui@fst.utm.tn}
\author[S. Jbari]{}
\email{sirine.jbari@fst.utm.tn}
\author[S. Tayachi]{}
\email{slim.tayachi@fst.rnu.tn}
\date{\today }
\subjclass[2020]{Primary: 35G20, 35G25, 35Q55. Secondary: 35Q70, 35Q35 }
\keywords{Inhomogeneous  Nonlinear Schr\"odinger  Equation, Local existence, Well-posedness,  Blow-up rate, Life span, Global existence, Decaying solution, Scattering theory.}
\begin{document}
\maketitle
\centerline{\scshape Lassaad Aloui} \medskip
{\footnotesize \centerline{Universit\'e de Tunis El Manar, Facult\'e des Sciences de Tunis,}  \centerline{ D\'epartement de
Math\'ematiques, Laboratoire  \'Equations aux} \centerline{ D\'eriv\'ees
Partielles LR03ES04,  2092 Tunis,
Tunisia}}\medskip
\centerline{\scshape Sirine Jbari} \medskip
{\footnotesize \centerline{Universit\'e de Tunis El Manar, Facult\'e des Sciences de Tunis,}  \centerline{ D\'epartement de
Math\'ematiques, Laboratoire  \'Equations aux} \centerline{ D\'eriv\'ees
Partielles LR03ES04,  2092 Tunis,
Tunisia}}

\medskip
\centerline{\scshape Slim Tayachi}\medskip
{\footnotesize \centerline{Universit\'e de Tunis El Manar, Facult\'e des Sciences de Tunis,}  \centerline{ D\'epartement de
Math\'ematiques, Laboratoire  \'Equations aux} \centerline{ D\'eriv\'ees
Partielles LR03ES04,  2092 Tunis,
Tunisia}}
\begin{abstract}
We are interested in the behavior of solutions to the damped  inhomogeneous  nonlinear Schr\"odinger  equation $ i\partial_tu+\Delta u+\mu|x|^{-b}|u|^{\alpha}u+iau=0$,  $\mu \in\C $, $b>0$, $a \in \C \  \mbox{such that}\   \Re \textit{e}(a) \geq 0$,  $\alpha>0$.    We establish lower and upper bound estimates of the life-span. In particular, for $a\geq 0$ we obtain explicit values  $a_*,\; a^*$ such that if $a<a_*$ then blow up occurs, while for $a>a^*,$ global existence holds. Also, we prove scattering results with precise  decay rates for large damping. Some of the results are new even for $b=0.$
\end{abstract}
 \section{Introduction}
 We consider the following inhomogeneous nonlinear damped Schr\"odinger equation
\begin{equation}
\label{INLS}
i\partial_tu+\Delta u+\mu|x|^{-b}|u|^{\alpha}u+iau=0, \quad
                 u(0,.)=u_0,
\end{equation}
where $u_0$ is a complex valued function, $ u_0\in H^s(\R^N)$, $s=0,1$,
$u=u(t,x)\in \C,\  t \in [0,\infty),\ \mu\in \C,\ x\in \R^N \setminus \{0\}$,  
$a \in \C \  \mbox{such that}\   \Re \textit{e}(a) \geq 0$
and
$$N\geq 1,\ N>2s   ,\  0<b<\min \{2,N-2s\}, 0< \alpha\leq \frac{4-2b}{N-2s}.$$
We interpret $0< \alpha< \infty$  when $s=1$ and $N=1,2$. The case $ \mu >0 $ (resp.
$ \mu <0 $) corresponds to the focusing (resp. defocusing ) case. We are interested in the impact of the damping term on  the global existence, the large time  behavior 
and the blow-up solutions for $\eqref{INLS}$. Our aim is to  extend the known results to $b>0$ and unify them with $b=0$. Some of our results are new even for $b=0$.

Equation \eqref{INLS} with $b=0$ and $a>0$ has various applications in the area of nonlinear optics, plasma
physics and fluid mechanics (see \cite{Fib,Kel,TW,TsutsumiM,ZSS}).
For $a=0$ and $b>0,$  \eqref{INLS} plays a crucial rule as a limiting case in the study of the nonlinear inhomogeneous Schr\"odinger equation
\begin{equation}
\label{INLSr}i\partial_t u +\Delta u  + K(x)|u|^\alpha u=0,
\end{equation}
with regular potential
$K(x)\sim |x|^{-b}$ as $|x|\to\infty$, which itself appears in diverse branches of physics such as in nonlinear optics \cite[Section 6]{Genoud}. As noted in \cite{Fib}, the damping (absorption) term plays an important effect in the physical model and it is better to not be neglected.  This motivated us to consider \eqref{INLS}. The difficulty in the study of \eqref{INLS} comes from the singularity of the potential at the origin. 

We are first interested in the local well posedness and the lower estimate of the life-span. Let us recall some known results for the case $a>0 $ and $b=0$. For this case, as mentioned in \cite{TsutsumiM,I,OT,Dinh2}, the local well posedness follows similarly as for the NLS equation. To our best knowledgement, no lower bound estimate for the life-span is known. On the other hand,   it is well known that  the global existence holds for the defocusing case and $0<\alpha<\frac{4}{N-2}$ or for the focusing case and $0<\alpha<\frac{4}{N}$. See for example \cite{I}.
For the case $ \frac{4}{N}\leq\alpha< \frac{4}{N-2}$, Ohta-Todorova \cite{OT} proved that for any $u_0\in H^1(\R^N) $ there exist $a^*(\|u_0\|_{H^1(\R^N)})>0$ such that  for all $ a>a^*(\|u_0\|_{H^1(\R^N)})$, the solution exists globally in time. See also \cite{SS,Fib} for numerical results. Recently, we have learned in \cite{Dinh2} that for the case $\alpha=\frac{4}{N-2}$ there exists $ a^*(u_0)>0$ such
that for all $ a>a^*(u_0)$  
the corresponding solution is global.

The main contribution of this paper is to develop the theory of the damped inhomogeneous nonlinear  Schr\"odinger equation \eqref{INLS}
 and establish a lower bound estimate of the life-span. To do this, we use the strategy of  \cite{AT} (see also \cite{TW1,TW2}) doing a fixed point argument in suitable Sobolev-Lorentz spaces. This allows us to obtain at once the life-span estimate and a lower bound of the damping term giving global existence.
 
We study  the problem \eqref{INLS} via its integral version
\begin{equation}
\label{intINLS}
u(t)=e^{t(i\Delta-a)}u_0+i\mu\int_0^te^{(t-\sigma)(i\Delta-a)}(|x|^{-b}|u|^\alpha u(\sigma))d\sigma.
\end{equation}
To give our result we need the following definition.
\begin{definition}
We say that a pair $ (q,r) $ is admissible if
$$ \frac{N}{2}=\frac{N}{r}+\frac{2}{q}, $$ 
and $ 2\leq r \leq \frac{2N}{N-2} $ ($ 2\leq r \leq \infty $ if $ N=1 $ and  $ 2\leq r<\infty  $ if $ N=2 $).
\end{definition}
We introduce the following particular admissible pair $ (\gamma,\rho)$ defined by
  \begin{equation}\label{admissible}
      \gamma=\frac{4(\alpha +2)}{\alpha (N-2s)+2b},\qquad  \qquad \rho=\frac{N(\alpha +2)}{N+\alpha s-b}
  \end{equation} 
  and
  \begin{equation}\label{theta3}
 \theta =\frac{4}{4-2b-\alpha(N-2s) }.
\end{equation}
We denote the Sobolev-Lorentz space by 
$$ W^{s,p}_q(\R^N):=  (I-\Delta)^{-\frac{s}{2}} L^{p,q}(\R^N),$$
$$ \dot{W}^{s,p}_q(\R^N):=(-\Delta)^{-\frac{s}{2}} L^{p,q}(\R^N), $$
 where $ L^{p,q}(\R^N)$ is the Lorentz space. We will recall the
definition of these spaces in Section 2. 

We have obtained the following.
\begin{theorem}[Lower-bound of the life-span]\label{ge}
Let $ N\geq 1$, $ \mu \in \C $, $s=0,1 $, $ s< \frac{N}{2}$  and $a\in \C $ such that $\Re e(a)> 0$. Assume that
$$ 0<b<\min \{2,N-2s\}\   \mbox{and}\   0< \alpha \leq \frac{4-2b}{N-2s}.  $$
Then  for every $u_0\in H^s(\R^N)$
 there exist   $T_{\max}(u_0)>0$  and a unique maximal solution $ u $ of the equation \eqref{intINLS} in  $ C([0,T_{\max}(u_0));H^s(\R^N))\cap L^{\gamma}_{loc}([0,T_{\max}(u_0));W^{s,\rho}_2(\R^N))$. Moreover, we have
the following properties.
  \begin{itemize}
 \item[(i)]
 $ u\in  L_{loc}^{p}([0,T_{\max}(u_0));W^{s,l}_2(\R^N)) $ for every admissible pair $ (p,l) $.
 \item[(ii)] (Blow-up criterion).
  If  $ T_{\max}(u_0)<\infty $,  then
 $$\|u\|_{L^{{\gamma}}([0, T_{\max}(u_0)):W^{s,\rho}_2(\R^N))}=\infty ,$$ 
 where $ (\gamma,\rho)  $ is given by (\ref{admissible}).
 \item[(iii)]
If $ \alpha <\frac{4-2b}{N-2s}$ and $T_{\max}(u_0)<\infty$, 
  then 
   $$ \lim_{t\to T_{\max}(u_0) }\|u(t)\|_{\dot{H}^s(\R^N)}=\infty .$$
\item[(iv)](Lower estimate of the life-span). For $0<\alpha<\frac{4-2b}{N-2s}$, there exists a constant $C=C(b,\alpha, N)>0 $ such that, 
$$T_{\max}(u_0)> \frac{1}{\alpha\theta \Re e(a)}\log \left(\frac{C\|u_0\|_{\dot{H}^s(\R^N)}^{\alpha\theta}}{\left(C\|u_0\|_{\dot{H}^s(\R^N)}^{\alpha \theta}-\alpha\theta\Re e(a)\right)_+}\right).$$ 
In particular, if 
    $$\Re e(a)\geq \displaystyle\frac{C}{\alpha\theta }\|u_0\|_{\dot{H}^s(\R^N)}^{\alpha\theta},$$
  then the solution  exists globally in time.
 \item[(v)](Critical case). Let $\alpha=\frac{4-2b}{N-2s}$. Then there exists $T_0(b,N,s,u_0)>0$ such that for all $\Re e(a)\geq 0$ we have  $T_{\max}(u_0)>T_0$. 
\item[(vi)] (Unconditional Uniqueness).  If $s=1$ and $N\geq 3$ or if $s=0 $ and
under one of this assumptions \\ 
 $\left\{
    \begin{aligned}
     &  N=1,\  0<b<\frac{1}{2},\  0<\alpha< 1-2b,\\
     &N=2,\  0<b<1,\  0<\alpha<1-b,\\
     &N\geq3,\  0<b<1,\  0<\alpha<\frac{2-2b}{N}. 
    \end{aligned}
  \right.$\\
  The equation has a unique solution in $ C([0,T_{\max}(u_0));H^s(\R^N) )$.
   \end{itemize}
\end{theorem}
\begin{remark}$\;${\rm
\begin{enumerate}
 \item[1)] The  result  of Theorem~\ref{ge} remains valid for the case $ b=0  $. In this
context, Theorem~\ref{ge} unifies the results for $b=0$ and $b>0$. 
\item[2)]
Our lower bound estimate of the life-span in Theorem~\ref{ge} (iv) is new even for $b=0$.
\item[3)] In Theorem~\ref{ge} (iv), the lower bound of $\Re e(a)$ giving global existence depends only on $\|u_0\|_{\dot{H}^s(\R^N)}$. In the case $b=0$,  a lower bound depending on $\|u_0\|_{H^s(\R^N)}$ is obtained in \cite{OT}. In  this context our result is new even for $b=0$.

\end{enumerate}

 }\end{remark}

 Our next result concerns with the local well-posedness   for $u_0\in H^1(\R^2)$.  
The author of  \cite{Guzman} has established the local well-posedness for  $0<b<\frac{2}{3}$. Recently, Dinh \cite{Dinh} improve  Guzm\'an's result on the rang of $b$ which can be extended to $0<b<1$.
 We revisit the local well-posedness in $ H^1(\R^2)$ and we obtain more precise regularity results, in addition we give a simple  proof for the case $a=0$. 
\begin{theorem}[Lower-bound of the life-span:   case $H^1(\mathbb{R}^2)$]\label{N2} Let $N=2$, $\mu \in \C$ and $a\in \C $ such that $\Re e(a)\geq 0$. Assume that
 $$0<b<1\  \mbox{and}\  0<\alpha<\infty.$$
 Then  for every $u_0\in H^1(\R^2)$
 there exist   $T_{\max}(u_0)>0$  and a unique maximal solution $ u $ of the equation \eqref{intINLS} in  $ C([0,T_{\max}(u_0)); H^1(\R^2))$. Moreover, we have
the following properties.
\begin{itemize}
 \item[(i)] $ u\in  L_{loc}^{p}([0,T_{\max}(u_0));W^{1,l}_2(\R^2)) $ for every admissible pair $ (p,l) $.
  \item[(ii)]If $T_{\max}(u_0)<\infty$, 
  then 
   $$ \lim_{t\to T_{\max}(u_0) }\|u(t)\|_{H^1(\R^2)}=\infty .$$
 \item[(iii)](Lower estimate of the life-span). For $0<\tau< \min\{ \frac{1-b}{\alpha+1}, \frac{\alpha}{\alpha+1}\},$ there exists a constant  $C>0$ such that 
 $$T_{\max}(u_0)
> \frac{2-b-(\alpha+1)\tau}{2\Re e(a) \alpha }\log\left(\frac{C\|u_0\|_{H^1(\R^2)}^{\frac{2\alpha}{2-b-(\alpha+1)\tau}}}{\left(C\|u_0\|_{H^1(\R^2)}^{\frac{2\alpha}{2-b-(\alpha+1)\tau}}-
\frac{2\alpha \Re e(a)}{2-b-(\alpha+1)\tau}\right)_+}\right),$$
 In particular, if $$ \Re e(a) \geq   \frac{2-b-(\alpha+1)\tau}{2\alpha}C\|u_0\|_{H^1(\R^2)}^{\frac{2\alpha}{2-b-(\alpha+1)\tau}},$$ then the solution  exists globally in time.
 \end{itemize}
\end{theorem}

Now, we turn to study the scattering problem for \eqref{INLS}. We give the following definition, see \cite{I}.
\begin{definition}
\label{def1}
Let $\psi : [0,\infty) \rightarrow [1,\infty)$ be nondecreasing for large time $t > 0$. We say that a global solution u to \eqref{INLS} with initial data $u_0\in H^s(\R^N),\; s=0,\; 1$, exponentially scatters with decay order $ \psi^{-1}$  if there exists $u_{+} \in H^s(\R^N)$ such that
\begin{equation}
    \lim_{t\rightarrow \infty}\psi(t)\e^{at} \|u(t)-\e^{t(i\Delta-a)}u_+\|_{H^s(\R^N)}=0,
\end{equation}
where $\e^{t(i\Delta-a)}$ is the semigroup generated by the maximal dissipative operator $i\Delta-a$.
\newline When $\psi=1$,  we say exponential scattering for
simplicity. That is we have 
\begin{equation}\label{expsca}
    \lim_{t\rightarrow \infty}\e^{at} \|u(t)-\e^{t(i\Delta-a)}u_+\|_{H^s(\R^N)}=0.
\end{equation}
\end{definition}

Regarding the scattering for \eqref{INLS} with $b=0$,  in \cite{I} it is shown that  for the defocusing case and $0<\alpha<\frac{4}{N-2}$ or for the focusing case and $0<\alpha<\frac{4}{N}$ all solutions, which are already global, scatter. See also \cite{OT,Dinh2} for the scattering with $b=0$ and large damping. We have obtained the following for $b>0.$

\begin{theorem}[Scattering for large damping]
\label{scatteringforlargea}
Let $\mu,\; a\in \mathbb{C}$ with $\Re e (a)>0$. Assume that
$$
N=2,\ s=1,\, 0<b<1\  \mbox{and}\  0<\alpha<\infty\ $$
\mbox{or} 
$$N \geq 1,\ s=0,\, 1,\, s<\frac{N}{2},\,  0<b<\min \{2,N-2s\}\   \mbox{and}\    0 < \alpha \leq \frac{4-2b}{N-2s}.
$$
Let $u_0\in H^s(\R^N)$ and $u$ be the unique maximal solution of the equation \eqref{intINLS}. Then the following hold.
\begin{itemize}
\item[(i)](Subcritical  case). Let $0<\alpha<\frac{4-2b}{N-2s}.$ Then there exists $a^*=a^*(\|u_0\|_{H^s},b,\alpha , N)\in (0,\infty)$ such that for $\Re e(a)\geqslant a^*,$ the solution  $u$ exists globally in time  and scatters exponentially.
\item[(ii)](Critical case). Let $\alpha=\frac{4-2b}{N-2s}.$ Then there exists $a^*=a^*(u_0,b,\alpha , N)\in (0,\infty)$ such that for $\Re e(a)\geqslant a^*,$ the solution  $u$ exists globally in time  and scatters exponentially.

\end{itemize}

\end{theorem}

We now will give other  results on the scattering for $\mu\in \R$ with precise estimates. For this we define the following quantities.
\begin{align}
    &\mbox{$L^2$-norme}:\  M(f):=\|f\|^2_{L^2(\R^N)}.\\
    &\mbox{Energy}:\ E(f):=\frac{1}{2}\|\nabla f\|_{L^2(\R^N)}^2-\frac{\mu}{\alpha +2}\||x|^{-b}|f|^{\alpha +2}\|_{L^1(\R^N)}.
\end{align}
However, it is easy to prove that if $u$ is a smooth solution of \eqref{INLS} on $[0,T)$, then for all $t\in [0,T)$ we have that
\begin{align}
& \partial_t M(u(t))=-2\Re \textit{e}(a)\|u(t)\|^2_{L^2(\R^N)},\label{masse}\\
&\partial_tE(u(t))=:-\Re \textit{e}(a)K(u(t)),\label{energy}
\end{align}
where
\begin{equation}
\label{energy1}
K(f)= \|\nabla f\|^2_{L^2(\R^N)}-\mu\||x|^{-b}|f|^{\alpha +2}\|_{L^1(\R^N)}.
\end{equation}

Exceptionally, 
the $L^2$-critical case, $\alpha=\frac{4-2b}{N}$,  $\mu=1,$
was treated by Genoud \cite{Genoud}. The author 
also showed how small should be the initial data to have global well-posed in $ H^1(\R^N)$ assuming $\|u_0\|_{L^2(\R^N)} <\|Q\|_{L^2(\R^N)}$ 
where $Q$ is the unique nonnegative, radially symmetric, decreasing solution of the ground state equation
\begin{equation}
\label{groundstate}\Delta Q-Q+|x|^{-b}|Q|^{\frac{4-
2b}{N}}Q=0.
\end{equation} 
 The next Theorem extends the well-known result of the behavior of the global solutions of  the classical ($b=0$) nonlinear damped Schr\"odinger equation.    For simplicity we will denote the norm in $L^p(\R^N),$ $\|\cdot\|_{L^p(\R^N)}$ by $\|\cdot\|_p.$
\begin{theorem}[Global existence and scattering]\label{sca1}
Let $ N\geq 3 $, $a\in\C$ $\mbox{such that}\ \Re e(a)>0 $, $u_0\in H^1(\R^N)$  and $ 0<b<\min\{2,N-2\} $.
\begin{itemize}
\item[1)]
Let  $\mu<0$ and   $ 0 <\alpha< \frac{4-2b}{N-2}$.
 Then, the solution of
\eqref{intINLS} exists globally in time and exponentially  scatters. In addition for all $t\geq 0$,
\begin{equation}\label{a}
    \|u(t)\|_2^2= \e^{-2\Re e(a)t}\|u_0\|_2^2,
\end{equation}
$$ E(u(t))\leq \e^{-2\Re e(a)t} E(u_0),$$
 $$\int_0^\infty \e^{2\Re\textit{e}(a)s} \||x|^{-b} |u(s)|^{\alpha +2}\|_1 ds\leq  \frac{\alpha+2}{\Re e(a)|\mu|\alpha} E(u_0).$$
\item [2)]Let $ \mu >0  $ and  $ 0 <\alpha< \frac{4-2b}{N}$. Then,   all  solutions of
\eqref{intINLS} exist globally in time, exponentially scatter and satisfy (\ref{a}).
In addition there exists a constant  $C=C(b,\alpha,a,N)>0$  such that  for all $t\geq 0$, 

\begin{equation}\label{b}
 \|\nabla u(t)\|_2^2
 \leq  C \e^{-2\Re\textit{e}(a)t}
 \left( \|\nabla u_0\|_2^2+\| u_0\|_2^{\beta}\right),
 \end{equation}
 \begin{align}\label{c}
   &\||x|^{-b} |u(t)|^{\alpha +2}\|_1
   \leq
   C\e^{-2\Re\textit{e}(a)t}\left(\|\nabla u_0\|_2^2+\| u_0\|_2^{\beta}+\||x|^{-b}|u_0|^{\alpha +2}\|_1\right),\\
   &\int_0^\infty \e^{2\Re\textit{e}(a)s} \||x|^{-b} |u(s)|^{\alpha +2}\|_1 ds
   \leq
   C\left(\|\nabla u_0\|_2^2+\| u_0\|_2^{\beta}+\||x|^{-b}|u_0|^{\alpha +2}\|_1\right), \label{d}
    \end{align}
    where
    $$\beta=\frac{2(4-2b-(N-2)\alpha)}{4-2b-N\alpha}.$$
\item[3)] Let $\mu>0$ and  $ \alpha =\frac{4-2b}{N} $. Assume that $\|u_0\|_2<\mu^{-1/\alpha}\|Q\|_2$ where $Q$ is the ground state of \eqref{groundstate}. Then  the global solution $u$  of \eqref{intINLS}  exponentially scatters and satisfies (\ref{a}). In addition, for all $t\geq0$, 
\begin{equation}
\label{intg1}
    \|\nabla u(t)\|_2^2\leq \e^{-2\Re\textit{e}(a)t} F(u_0),
    \end{equation}
\begin{equation}
\label{intg2}
 \e^{2\Re\textit{e}(a)t} \||x|^{-b} |u(t)|^{\alpha +2}\|_1+\Re\textit{e}(a)\alpha\int_0^\infty  \e^{2\Re\textit{e}(a)s} \||x|^{-b} |u(s)|^{\alpha +2}\|_1ds\leq CF(u_0),
    \end{equation}
 where 
 $$F(u_0)={2E(0)\| Q\|_2^\alpha \over \| Q\|_2^\alpha-\mu \| u_0\|_2^\alpha}\exp \left({\mu\|u_0\|_2^{\alpha}\over \| Q\|_2^\alpha-\mu \| u_0\|_2^\alpha}  \right),\;  C={N+2-b\over \mu N}.$$
\end{itemize}
\end{theorem}
\begin{remark}
{\rm For  $\mu>0$ and $ \alpha ={(4-2b)}/{N},$ it is not clear that  any global solution exponentially scatters even for $b=0$.}
\end{remark}

We have the following result about de decay rate for the solutions that scatter.
\begin{theorem}[Decay rate of the scattering solutions]
\label{decayratescattering}
Let $ N\geq 3, $ $\mu \in \C, \   0\leq b<\min\{2,N-2\}$ and  $$0<\alpha\leq  \frac{4-2b}{N-2} .$$ If the global solution
of \eqref{equintegrale} scatters in $H^1(\R^N)$, then
$$ \lim_{t\rightarrow \infty} \e^{\alpha a t}\|v(t)-\e^{it\Delta}v_+\|_{H^1(\R^N)}=0.$$ 
\end{theorem}
\begin{remark} {\rm The previous Theorem, which holds also for $b=0$, shows that $u$ exponentially scatters with decay order $\e^{-\alpha a t},$ according to Definition \ref{def1}. This solves the first part of the open problem 2 in \cite{I}}.
  \end{remark}

We now turn to the blowup results of this paper.  We begin by  recalling the known ones for the case  $a=b=0$. The first blowup result for the undamped NLS is due to  \cite{Glassey}.  We refer also the reader to the book  \cite{Cazenave} and some references therein for other results. In particular, it is known that  for  $ \frac{4}{N}\leq \alpha< \frac{4}{N-2},$ and  $$u_0\in \Sigma:=\{u\in H^1(\R^N);\; \||\cdot|u\|_{L^2(\R^N)}<\infty\}$$ when  one of the  following conditions holds:
\begin{itemize}
    \item[(i)] $E(u_0)<0$,
    \item[(ii)]  $E(u_0) =0$ and $ V_0< 0$,
    \item[(iii)]  $E(u_0) >0$ and $ V_0<-\sqrt{2E(u_0)I_0} $,
\end{itemize}
where $$V_0:=\Im m\int_{\R^N} x.\nabla u_0(x) \overline{u_0}(x)dx\; \mbox{and}\;  I_0:=\||\cdot|u_0\|_{L^2(\R^N)}^2.$$

 For the homogenous damped NLS, that is $a>0$ and $b=0,$ Tsutsumi \cite{TsutsumiM} proved that if    $ \frac{4}{N}<\alpha<\frac{4}{N-2},\, u_0\not\equiv 0$ and  $u_0\in \Sigma$ satisfies
$ E(u_0)\leq 0$ and  $V_0< 0$
then the corresponding solution blows up in finite time for $$a< \frac{(N\alpha-4)|V_0|}{\alpha I_0},$$ and we have  $$T_{\max}(u_0)\leq {N\alpha-4\over 4\alpha a}\log\left[{(N\alpha-4)V_0\over (N\alpha-4)V_0+\alpha a  I_0}\right].$$

 In \cite{OT} Ohta and Todorova  improved Tsutsumi’s result  by showing that if $ \frac{4}{N}<\alpha<\frac{4}{N-2},\, u_0\not\equiv 0$ and $u_0\in \Sigma $ satisfies one of the above conditions (i)-(iii)
then there exists $a_*(u_0)>0$ such that for all $0\leq a<a_*(u_0) $, the corresponding solution blows up
in finite time.
Furthermore, 
similar results  are proved in the mass-critical case  for $N\geq 1$ and in the energy-critical for $N\geq 3$   by Dinh \cite{Dinh2}.
Our last aim is to prove a blow up result.

In this paper we are interested  in the case  $a>0,\; b>0.$ To our best knowledge, no  blowup results are known for  this case.  Recall that in the defocusing case with $0<\alpha <\frac{4-2b}{N-2} $  and  for the focusing case with $0<\alpha <\frac{4-2b}{N} $ , we have $ T_{\max}(u_0)=\infty $ for all $u_0 \in H^1  $. Therefore, we restrict our attention to the case $ \alpha > \frac{4-2b}{N}  $ and $ \mu =1$. For simplicity, we put
 We prove the following.

\begin{theorem}[Upper bounds of the  life-span]\label{blow}
Let $\mu=1,$ $a>0$ and 
$$
N=2,\  0<b<1\  \mbox{and}\  2-b<\alpha<\infty\ $$
\mbox{or} 
$$N \geq 3,\  0<b<\min \{2,N-2\}\   \mbox{and}\    \frac{4-2b}{N} < \alpha < \frac{4-2b}{N-2}.
$$
  Let $u_0 \in \Sigma$ and  $T_{\max}(u_0)$ be the existence time of the maximal solution  to \eqref{intINLS} with initial data $u_0$. Then the following hold.
  \begin{itemize}
\item[(i)] If $E_0<0,$ $\gamma\in \left[0,{\sqrt{V^2_0-4E_0I_0}-V_0\over 2I_0}\right]$ then $$T_{\max}(u_0)\leq {2I_0\over \sqrt{16(V_0^2-2I_0E_0)-I_0^2\gamma^2}\;-I_0\gamma-4V_0},$$
where
$$E(u_0)=E_0\; \mbox{and}\;\gamma=\frac{4\alpha}{N\alpha -4+2b}a.$$
\item[(ii)] If $E_0=0 $, $V_0<0 $, and $\gamma\in \left[0, {4|V_0|\over I_0}\right)$ then $$T_{\max}(u_0)\leq {1\over \gamma}\log\left(1+{I_0\gamma\over |I_0\gamma+4V_0|}\right), $$
\item[(iii)] If $E_0>0 ,$  $-2\sqrt{E_0I_0}\leq V_0<-\sqrt{2E_0I_0} $ and $\gamma\in \left[0,\displaystyle {4\over I_0|V_0|}\left(V_0^2-2E_0I_0\right)\right)$ then $$T_{\max}(u_0)\leq {1\over \gamma}\log\left({16E_0\over I_0\gamma^2+4V_0\gamma+16E_0}\right).$$
\item[(iv)] If $E_0>0 ,$  $V_0<-2\sqrt{E_0I_0} $ and  $\gamma\in \left[0,\displaystyle {2\over  I_0}\left(|V_0|-{\sqrt{V^2_0-4E_0I_0}}\right)\right)$ then $$T_{\max}(u_0)\leq {1\over \gamma}\log\left({16E_0\over I_0\gamma^2+4V_0\gamma+16E_0}\right).$$
\item[(v)] If $E_0>0 ,$  $V_0<-2\sqrt{E_0I_0} $ and  $\gamma\in \left[\displaystyle {2\over  I_0}\left(|V_0|-{\sqrt{V^2_0-4E_0I_0}}\right),\displaystyle {2\over  I_0}\left(|V_0|+{\sqrt{V^2_0-4E_0I_0}}\right)\right)$ then $$T_{\max}(u_0)\leq {1\over \gamma}\log\left(1+{ I_0\gamma^2\over |I_0\gamma^2+4V_0\gamma+16E_0|}\right).$$
\end{itemize}   
\end{theorem}
\begin{remark}
{\rm Theorem \ref{blow} shows  that more the value of the variance becomes negative,  wider is the range of values of the damping giving blow up.}
\end{remark}

Let us give some examples of initial data $u_0\in \Sigma$  such that one of the conditions (i)-(iv)  holds. 
For the first case it is enough to take $u_0(x)=\lambda \phi(x)$ such that  $\phi \in \Sigma(\R^N)$ with  $\lambda>0$ sufficiently large. For the other cases we may take 
 $u_0(x)=\sigma \psi(x)$ and  $\psi(x)=\e^{-i|x|^2} \Theta (x)$ such that $\Theta \in C_0^\infty(\R^N)$, $ \Theta \neq0$ and $\sigma>0 $ satisfies
$$\sigma^\alpha=\frac{\alpha +2}{2}\frac{\|\nabla \psi\|_{L^2(\R^N)}^2}{\||\psi|^{\alpha +2}|x|^{-b}\|_{L^1(\R^N)}},$$ 
for the second case and 
$$(\alpha+2)\frac{\left(\frac{1}{2}\|\nabla \psi\|_{L^2(\R^N)}^2 -4\|x\psi\|_{L^2(\R^N)}^2\right)}{\| |\psi|^{\alpha+2}|x|^{-b}\|_{L^1(\R^N)}} <\sigma^{\alpha}<\frac{\alpha+2}{2}\frac{\|\nabla \psi\|_{L^2(\R^N)}^2}{\| |\psi|^{\alpha+2}|x|^{-b}\|_{L^1(\R^N)}}, 
$$
for the third case.
\par The rest of this paper is organized as follows. In the next section, we introduce some notation and recall Strichartz
estimates in the Lorentz space  for the linear Schr\"odinger equation. In Section 3, we prove the local well-posedness and global existence. Section 4, 
is devoted to the scattering results, that is  the proof of Theorems ~\ref{scatteringforlargea}, ~\ref{sca1} and \ref{decayratescattering}. We will give the proof of Theorem~\ref{blow} in Section 5.
In the rest of the paper  we  abbreviate $ L^{r,2}(\R^N)$, $ W^{1,r}_2(\R^N)$ and  $H^1(\R^N) $ to  
 $L^{r,2}$,  $ W^{1,r}_2 $ and   $H^1 $, respectively. The notation $A\lesssim B$ for positive numbers $A$ and $B$, means that there exists a positive constant $C$ such that $A\leq CB.$ If $A\lesssim B$ and $B\lesssim A,$ we write $A\sim B.$ We denote $a_+:=\max(0,a)$ for  a real number $a$. Also, $C$ will denotes a constant which may be different at different places. 

\section{Preliminaries}
\par In this section we recall the definitions of the Lorentz and the Sobolev-Lorentz  spaces. We also recall the Strichartz estimates  in the Lorentz spaces. Let $ f $ be a measurable function on $\R^{N}$.  
We define  on $[0,\infty) $ the distribution function of $ f $, which we denote $f_*$, by 
$$ f_*(\lambda):= \mu (\{x\in \R^N; |f(x)|> \lambda \}),\;  \lambda \geq 0, $$ 
 where $ \mu$ represents the  Lebesgue measure. We define the rearrangement $ f^*:[0,\infty)\rightarrow [0,\infty] $  of measurable function  $ f $ by 
 $$ f^* (t):=\inf \{\lambda>0; f_*(\lambda)\leq t\}, \;  \ t\geq 0.$$
 Moreover, $ f^{**}:(0,\infty)\rightarrow [0,\infty] $ denotes the average function of $ f^* :$ 
 $$f^{**}(t):=\frac{1}{t}\int_0^t f^*(s)ds, \;    t>0.$$
For $ 1\leq p<\infty $ and $ 1\leq q \leq \infty $, the Lorentz space $ L^{p,q}(\R^N) $ is defined as the set of  measurable functions $f$  such that the following  quasi-norm is finite:
\begin{equation}\label{norme}
\|f\|_{L^{p,q}}:=\left\{
    \begin{aligned}
     &  \left( \int_0^\infty (t^{\frac{1}{p}}f^*(t))^q \frac{dt}{t} \right)^\frac{1}{q}\ \  \mbox{if}\qquad   1\leq q<\infty ,\\
     &\sup_{t>0}(t^{\frac{1}{p}}f^*(t))\ \  \mbox{if} \qquad   q=\infty.
    \end{aligned}
  \right.
\end{equation}
If we replace$ f^* $ by $ f^{**} $ in (\ref{norme})  we get  an equivalent norm for $ L^{p,q}(\R^N) ,$ $ p> 1 $ and in this case, $L^{p,q}(\R^N) $ is a Banach space. The Sobolev-Lorentz spaces  are defined by 
   $$ W^{1,p}_q(\R^N):=  (I-\Delta)^{-\frac{1}{2}} L^{p,q}(\R^N),\; \dot{W}^{1,p}_q(\R^N):=  (-\Delta)^{-\frac{1}{2}} L^{p,q}(\R^N)$$
   with $ N>1 $ , $ 1\leq p<\infty $
and $ 1\leq q \leq \infty $,   
equipped with the norms $ \|f\|_{W^{1,p}_q(\R^N)} := \|(I-\Delta)^{\frac{1}{2}}f\|_{L^{p,q}(\R^N)}$, $ \|f\|_{ \dot{W}^{1,p}_q(\R^N)} := \|(\Delta)^{\frac{1}{2}}f\|_{L^{p,q}(\R^N)}$. 
If $ 1<p<\infty $ and $ 1\leq q_1<q_2\leq \infty $, then 
$$\|f\|_{L^{p,q_2}(\R^N)}\lesssim \|f\|_{L^{p,q_1}(\R^N)}. $$ We  now recall the H\"older inequality in the Lorentz spaces. If $ h=fg $
with $ f \in L^{p_1,q_1}(\R^N) $ and $  g\in L^{p_2,q_2}(\R^N)  $ where 
$$ \frac{1}{p_1}+\frac{1}{p_2}<1,
 $$
  then $ h\in   L^{r,s}(\R^N) $ where 
  $$ \frac{1}{p_1}+\frac{1}{p_2}=\frac{1}{r}, $$
and $ s\geq 1 $ is any number such that 
$$ \frac{1}{q_1}+\frac{1}{q_2}\geq \frac{1}{s}. $$   
Moreover $$\|h\|_{L^{r,s}(\R^N)}\lesssim \|f\|_{L^{p_1,q_1}(\R^N)} \|g\|_{L^{p_2,q_2}(\R^N)}.$$
If $ f\in L^{p,q_1}(\R^N) $ and $  g\in L^{p',q_2}(\R^N)  $ where
$$  \frac{1}{p}+\frac{1}{p'}=1,$$
$$ \frac{1}{q_1}+ \frac{1}{q_2}\geq 1, $$
then $ h\in L^1(\R^N) $ and  
$$ \|h\|_{L^1(\R^N)}\leq \|f\|_{L^{p,q_1}}\|g\|_{L^{p',q_2}(\R^N)}. $$

We recall the homogeneous Sobolev-Lorentz embedding:
 $\dot{W}^{1,p}_q(\R^N)\hookrightarrow L^{\tilde{p},q}(\R^N)$, where $1<p<N $,  $1\leq q\leq \infty $, and 
 $$ \frac{1}{\tilde{p}}=\frac{1}{p}-\frac{1}{N}.$$
 
 By the well known Sobolev embedding $H^s(\R^N)\hookrightarrow L^p(\R^N) $ and interpolation, we have the following
 $$ H^s(\R^N)\hookrightarrow L^{p,2}(\R^N),\qquad s\geq 0,\qquad \frac{1}{2}-\frac{s}{N}\leq \frac{1}{p}\leq \frac{1}{2},\qquad p<\infty. $$
 
 \begin{proposition}[Gagliardo Nirenberg inequality (\cite{HYZ})]
  Let  $1\leq p,\; q,\; p_2,\; q_1,\; q_2<\infty $, $0<\beta<q $, $0<s<N $ and $1<p_1<\frac{N}{s}$. Then there exists a constant $A>0 $ such that 
  $$\|u\|_{L^{p,q}(\R^N)}\leq A \|(-\Delta)^{\frac{s}{2}}u\|^{\frac{\beta}{q}}_{L^{p_1,q_1}(\R^N)}\|u\|_{L^{p_2,q_2}(\R^N)}^{\frac{q-\beta}{q}} $$
  holds for 
  \begin{align*}
      &\frac{\beta}{q_1}
+\frac{q-\beta}{q_2}=1,\\
&\beta \left(\frac{1}{p_1}-\frac{s}{N} \right)+(q-\beta)\frac{1}{p_2}=\frac{q}{p}.
\end{align*}
 \end{proposition}

We have the following Strichartz estimates in Lorentz spaces.
\begin{proposition}[Strichartz's estimates (\cite{KT})] Let $ (q,r) $ and $ (\tilde{q},\tilde{r}) $ be admissible pairs and let I be an interval of $ \R $.  Then there exist two constants $ c_1 $ and $ c_2 $  such that we have the following estimates 
 $$\|\e^{it\Delta }\varphi\|_{L^q(I,L^{r,2}(\R^N))} \leq c_1 \|\varphi\|_{L^2(\R^N)}, \ \ \mbox{for every}\ \varphi \in L^2(\R^N).  $$
  $$\left\|\int_0^t\e^{i(t-s)\Delta }f\right\|_{L^q(I,L^{r,2}(\R^N))} \leq c_2 \|f\|_{L^{\tilde{q}'}(I,L^{\tilde{r}',2}(\R^N))}, \ \ \mbox{for every}\ f \in L^{\tilde{q}'}(I,L^{\tilde{r}',2}(\R^N)) .$$
\end{proposition}

We end this section by introducing  the following transformation :
\begin{equation}\label{transformation}
      v(t,x):=e^{at}u(t,x),
      \end{equation}
      which will be used in our proofs. It follows that  $ u $ is a solution  of  \eqref{INLS} if and only if  $ v $
 is a solution  of
 \begin{equation}\label{2}
   i\partial_tv+\Delta v+\mu h(t)|x|^{-b}|v|^{\alpha}v=0, \quad v(0)= v_0=u_0 ,
\end{equation}
where
$h(t)=\e^{-\alpha at}$.  
We study the previous equation via its Duhamel formulation
\begin{equation}\label{equintegrale}
 v(t)=\e^{it\Delta}u_0+i\mu \int_0^t \e^{i(t-s)\Delta}[h(s)|x|^{-b}|v(s)|^{\alpha}v(s)]ds.
  \end{equation}

\section{Local well-posedness and Global existence}
The purpose of this section is to study the local well-posedness of the Cauchy problem \eqref{INLS}  in $H^s(\R^N)$ for $s=0,1$ and $s<\frac{N}{2}$. The main tools are the Strichartz estimates in
Lorentz spaces. To do this we separate the cases $N\geq 3$ and $N=2$.
\subsection{ Case $N\geq 3$ }

 In order to prepare the proof of Theorem \ref{ge}, we establish the following estimates. 
 \begin{lemma}
 \label{estimationstri}
 Let $N\geq 1$, $s=0,1$, $ s<\frac{N}{2}$, $0<b<\min\{2,N-2s\}$, $0<\alpha \leq \frac{4-2b}{N-2s}$ and $ a\in \C$ such that $\Re\textit{e}(a)\geq0$.  Let $T>0,$ $(\gamma,\rho)$ be given by \eqref{admissible} and  $u,w\in L^\gamma(0,T,W^{s,\rho}_2)$. Then
              \begin{equation}\label{j1}
     \|h|x|^{-b}|u|^{\alpha}w\|_{L^{\gamma'}((0,T]:L^{\rho',2})}\lesssim \|h\|_{L^{\theta}([0,T])}\|(-\Delta)^{\frac{s}{2}} u\|^{\alpha}_{L^\gamma([0,T]:L^{\rho,2})}
\|w\|_{L^\gamma([0,T]:L^{\rho,2})}
 \end{equation}
  and  
     \begin{equation}\label{j2}
     \|h|x|^{-b-s}|u|^{\alpha}w\|_{L^{\gamma'}([0,T)]:L^{\rho',2})} \lesssim  \|h\|_{L^{\theta}([0,T))}\|(-\Delta)^{\frac{s}{2}} u\|^\alpha_{L^{\gamma}([0,T]:L^{\rho,2})}\|(-\Delta)^{\frac{s}{2}} w\|_{L^{\gamma}([0,T]:L^{\rho,2})},
     \end{equation}
     where $ \theta $ is given by \eqref{theta3}.
 \end{lemma}
 \begin{proof}
 We first prove \eqref{j1}. By the  H\"older inequality in Lorentz spaces we have
\begin{align*}
 \|h|x|^{-b}|u|^{\alpha}w\|_{L^{\gamma'}((0,T]:L^{\rho',2})}&\lesssim \||x|^{-b}\|_{L^{\frac{N}{b},\infty}}\|h|u|^{\alpha}w\|_{L^{\gamma'}([0,T]:L^{\eta,2})}\\
 &\lesssim \|h\|_{L^{\theta}([0,T])} \|u\|^\alpha _{L^{\gamma}([0,T]:L^{n_1,2})}\|w\|_{L^{\gamma}([0,T]:L^{\rho,2})},
\end{align*}
provided that
$$ \frac{1}{\eta}+\frac{b}{N}=\frac{1}{\rho'},\ 
     \frac{\alpha}{n_1}+\frac{1}{\rho}=\frac{1}{\eta},\ 
   \frac{1}{\rho}-\frac{s}{N} =\frac{1}{n_1},\ 
  \frac{\alpha}{\gamma}+\frac{1}{\gamma}+\frac{1}{\theta}=\frac{1}{\gamma'}.    
 $$

Since $0<s<\frac{N}{2}$ and $0<b<\min \{2,N-2s\}$, we can use the homogenous Sobolev embedding $\dot{W}^{s,\rho}_2  \hookrightarrow L^{n_1,2}$ to get
$$ \|h|x|^{-b}|u|^{\alpha}w\|_{L^{\gamma'}([0,T]:L^{\rho',2})}\lesssim \|h\|_{L^{\theta}([0,T])}\|(-\Delta)^{\frac{
s}{2}} u\|^\alpha_{L^{\gamma}([0,T]):L^{\rho,2})} \|w\|_{L^{\gamma}([0,T]:L^{\rho,2})}.$$
We now prove \eqref{j2}. Using H\"older's inequality  we have
\begin{align*}
 \|h|x|^{-b-s}|u|^{\alpha}w\|_{L^{\gamma'}([0,T):L^{\rho',2})}&\leq \||x|^{-b-s}\|_{L^{\frac{N}{b+s},\infty}}\|h|u|^\alpha w\|_{L^{\gamma'}([0,T):L^{\eta_2,2})}\\
&\lesssim \|h\|_{L^{ \theta}([0,T])}\|u\|^\alpha_{L^{\gamma}([0,T]:L^{n_1,2})}\|w\|_{L^{\gamma}([0,T]:L^{n_1,2})},
\end{align*}
provided that
$$
     \frac{1}{\eta_2}+\frac{b+s}{N}=\frac{1}{\rho'},\ 
     \frac{\alpha}{n_1}+\frac{1}{n_1}=\frac{1}{\eta_2},\ 
     \frac{\alpha}{\gamma}+\frac{1}{\gamma}+\frac{1}{\theta}=\frac{1}{\gamma'}.
     $$
Similarly as above, using the homogenous Sobolev embedding, we obtain
\begin{align*}
 \|h|x|^{-b-s}|u|^{\alpha}w\|_{L^{\gamma'}([0,T):L^{\rho',2})}
&\lesssim \|h\|_{L^{\theta}([0,T])}\|(-\Delta)^{\frac{s}{2}} u\|^{\alpha}_{L^{\gamma}([0,T]:L^{\rho,2})}\|(-\Delta)^{\frac{s}{2}} w\|_{L^{\gamma}([0,T]:L^{\rho,2})}.
\end{align*}

 \end{proof}
 \begin{proof}[Proof of Theorem~\ref{ge}] We use the  Banach fixed point theorem. 
Let $ M, T > 0$ to be chosen later and $(\gamma, \rho)$ be the admissible pair given by (\ref{admissible}). Consider the set
$$ E=\{v\in L^\gamma([0,T]:W^{s,\rho}_2), \|v\|_{L^\gamma([0,T]:\dot{W}^{s,\rho}_2)}\leq M \} $$
equipped with the distance
$$d(v,w)=\|v-w\|_{L^\gamma ([0,T]:L^{\rho,2})}.$$
Note that  $(E,d)$  is a complete metric space (see  \cite{AT}).
We Construct the local solution separately for the subcritical and critical cases.
 
\textbf{The subcritical case: $ \alpha<(4-2b)/(N-2s)$.}
Let us show that the mapping $\phi$ defined by
\begin{equation}\label{integrale}
 \phi(v)(t)=\e^{it\Delta}u_0+i\mu \int_0^t \e^{i(t-s)\Delta}[h(s)|x|^{-b}|v(s)|^{\alpha}v(s)]ds
  \end{equation}
  applies the space $ E$ into itself and it is a strict contraction.  Let $v\in E $. It follows from the Strichartz estimates that for every admissible pair $(p,l)$, we have
\begin{align*}
\|\phi (v)\|_{L^p([0,T]:\dot{W}^{s,l}_2)}&\leq \|\e^{it\Delta}u_0\|_{L^p([0,T]:\dot{W}^{s,l}_2)}+|\mu|\left\|\int_0^t \e^{i(t-s)\Delta}[h(s)|x|^{-b}|v(s)|^{\alpha}v(s)] ds \right\|_{L^p([0,T]:\dot{W}^{s,l}_2)}\\
&\lesssim\| u_0\|_{\dot{H}^s}+\|h|x|^{-b-s}|v|^{\alpha}v\|_{L^{\gamma'}([0,T]:L^{\rho',2})}+\|h|x|^{-b}|v|^{\alpha}(-\Delta)^{\frac{s}{2}} v\|_{L^{\gamma'}([0,T]:L^{\rho',2})}.
\end{align*}
Using (\ref{j1})  and   (\ref{j2}), we get
$$ \|\phi (v)\|_{L^p([0,T]:\dot{W}^{s,l}_2)}\lesssim \|u_0\|_{\dot{H}^s}+ \|h\|_{L^{\theta}([0,T])} \|v\|_{L^{\gamma}([0,T]:\dot{ W}_2^{s,\rho})}^{\alpha+1}.$$
If particular, if we choose $(p,l)=(\gamma,\rho) $,
then  $\phi(v)\in L^\gamma ((0,T];W^{s,\rho}_2) $ and we obtain the  following estimate
$$   \|\phi (v)\|_{L^\gamma([0,T];\dot{W}^{s,\rho}_2)}\lesssim \|u_0\|_{\dot{H}^s}+  M^{\alpha+1} \|h\|_{L^{\theta}([0,T])}.$$
Let $v,w\in E$. By the well known inequality
$$||v|^\alpha v-|w|^\alpha w|\lesssim (|v|^\alpha +|w|^\alpha) |v-w|,$$
and using the Strichartz estimates and the H\"older inequality, we obtain
\begin{align}
d(\phi (v),\phi (w))&\lesssim \||x|^{-b}h|v-w||v|^{\alpha}\|_{L^{\gamma'}([0,T]:L^{\rho',2})}+\||x|^{-b}h|v-w|||w|^{\alpha}\|_{L^{\gamma'}([0,T]:L^{\rho',2})}.\end{align}
Applying (\ref{j1}), we get 
$$d(\phi (v),\phi (w))\lesssim M^{\alpha}\|h\|_{L^{\theta}([0,T])} d(v,w).$$
Thus, there exists $K>0$ independent of $a$ such that 
$$
\left\{
    \begin{array}{l}
      \|\phi (v)\|_{L^\gamma([0,T]:\dot{W}^{s,\rho}_2)}\leq K\|u_0\|_{\dot{H}^s}+KM^{\alpha+1} \|h\|_{L^{\theta}([0,T])}, \\
     d(\phi (v),\phi (w))\leq KM^{\alpha}\|h\|_{L^{\theta}([0,T])} d(v,w).
    \end{array}
  \right.
 $$
 Finally, we need to choose $ M>0 $ and $ T > 0  $ such that
$$
\left\{
    \begin{array}{l}
     K\|u_0\|_{\dot{H}^s}+KM^{\alpha+1} \|h\|_{L^{\theta}([0,T])}\leq M, \\
      KM^{\alpha}\|h\|_{L^{\theta}([0,T])}<1 .
    \end{array}
  \right.
 $$
 For this we consider $M=2K\|u_0\|_{\dot{H}^s}$ and $T$ be such that 
 \begin{equation}\label{N4}
2^{\alpha}K^{\alpha+1} \|u_0\|_{\dot{H}^s}^{\alpha}\|h\|_{L^{\theta}([0,T])}\leq \frac{1}{2}.
 \end{equation}
 So the fixed point argument on $E$ with the above values of $ M$ and $T$ works and gives  a unique solution $v\in E$ of (\ref{2}).  
To show that $ v\in C([0,T];H^s)$ it suffices to prove that
$$  w(t)=\int_0^t \e^{-i\sigma\Delta}[h(\sigma)|x|^{-b}|v(\sigma)|^\alpha v(\sigma)]d\sigma \in C([0,T];H^s). $$
Let $s<\frac{N}{2} $ and $0\leq t<t'<T$. By Strichartz estimates and $(\ref{j1} )$, we get
 \begin{align*}
\|w(t)-w(t')\|_{H^s}&=\left\|\int_t^{t'} \e^{-i\sigma\Delta}[h(\sigma)|x|^{-b}|v(\sigma)|^\alpha v(\sigma)]d\sigma\right\|_{H^s}\\
&\lesssim \|h\|_{L^\theta((t,t'))} \|v\|_{L^{\gamma}((t,t'):W_2^{s,\rho})}.
\end{align*} 
 Since  $\|v\|_{L^{\gamma}((0,T):W_2^{s,\rho})} < \infty$, then if $ t \rightarrow t' $, $\|w(t)-w(t')\|_{H^s}\rightarrow 0.$ \\
By well known argument we may define the maximal time of existence $T_{\max}(u_0)>0$ and the maximal solution  of \eqref{intINLS} in  $ C([0,T_{\max}(u_0));H^s(\R^N))\cap L^{\gamma}_{loc}([0,T_{\max}(u_0));W^{s,\rho}_2(\R^N))$.

\textbf{The critical case : $ \alpha =\frac{4-
2b}{N-2s} $}.
Let $v\in E $. We  argue as above. Since $\theta=\infty$, then $\|h\|_{L^\theta([0,T])}=1$ and so we get 
 \begin{equation}
 \|\phi(v)\|_{L^\gamma([0,T]:\dot{W}^{s,\rho}_2)}\leq \|\e^{it\Delta}u_0\|_{L^\gamma([0,T]:\dot{W}^{s,\rho}_2)}+K  M^{\alpha +1}.
\end{equation}  
 Let $ v,w\in E $. Similarly,  we obtain 
    \begin{align*}
    d(\phi(v),\phi (w))\leq K M^{\alpha }d(v,w). 
    \end{align*}
    To apply fixed point argument we need that 
       $$
\left\{
    \begin{array}{l}
     \|\e^{it\Delta}u_0\|_{L^\gamma([0,T]:\dot{W}^{s,\rho}_2)}+K M^{\alpha +1}\leq M, \\
      K M^{\alpha }< 1.
    \end{array}
  \right.
 $$
 Since $\|\e^{it\Delta}u_0\|_{L^\gamma([0,T]:\dot{W}^{s,\rho}_2)}<\infty  $, then $ \|\e^{it\Delta}u_0\|_{L^\gamma([0,T]:\dot{W}^{s,\rho}_2)}\rightarrow 0 $ as $ T\rightarrow 0 $. Hence for all $ \varepsilon >0 $, there exists  $ T >0 $ such that $ \|\e^{it\Delta}u_0\|_{L^\gamma([0,T]:\dot{W}^{s,\rho}_2)}\leq \varepsilon $. Fix $ M=\frac{1}{(2K)^{\frac{1}{\alpha}}}$. For $ \varepsilon = \frac{M}{2}$, there exists $T>0$ such that $ \|\e^{it\Delta}u_0\|_{L^\gamma([0,T]:\dot{W}^{s,\rho}_2)}\leq \frac{M}{2}$. With such  a choice of $M$ and $T$ the above conditions are satisfied.  Consequently, we deduce the existence of a unique solution $ v\in E $ by the fixed point argument.

The proof of the blow up alternative, the regularity  and the unconditional uniqueness can be done as in \cite{AT}, so we omit it.

$(iv)$ We note that
\begin{equation*}
    \|h\|^\theta_{L^\theta([0,T])}=\frac{1}{\alpha\Re e(a)\theta}[1-\e^{-\alpha \Re e(a)\theta T}].
\end{equation*}
A simple calculation shows that if $\Re e (a)< \frac{(2^{\alpha+1}K^{\alpha+1}\|u_0\|_{\dot{H}^s}^\alpha)^{\theta}}{\theta\alpha} $, then 
\begin{equation*}
    T_{\max}(u_0)> \frac{1}{\alpha\theta \Re e(a)}\log \left(\frac{(2^{\alpha+1}K^{\alpha+1}\|u_0\|_{\dot{H}^s}^{\alpha})^{\theta}}{(2^{\alpha+1}K^{\alpha+1}\|u_0\|_{\dot{H}^s}^{\alpha})^{\theta}-\alpha\theta\Re e(a)}\right).
\end{equation*}
On the other hand, if $\Re e (a)\geq \frac{(2^{\alpha+1}K^{\alpha+1}\|u_0\|_{\dot{H}^s}^\alpha)^{\theta}}{\theta\alpha} $, then $T_{\max}(u_0)=\infty$.

$(v)$ It is clear from the construction of the local solution  for $\alpha=\frac{4-2b}{N-2}$ that  there exists $T_0>0$ independent of $a$ such that $T_{\max}(u_0) >T_0$. 
\end{proof}

\subsection{Case $N=2$ }
In this sub-section we prove the local existence of solution to equation \eqref{intINLS} in $ H^1(\R^2)$, that is Theorem~\ref{N2}.
\begin{proof}[Proof of Theorem~\ref{N2}]
Let $T,M > 0 $ to be chosen later and consider the set 
$$ E_2=\{v\in L^\infty([0,T];H^1),\  \|v\|_{L^\infty ([0,T];H^1)} \leq M \} $$
equipped with the distance
$$d(v,w)=\|v-w\|_{L^\infty ([0,T];L^2)}.$$
Applying Strichartz's estimate on \eqref{integrale}, we have
\begin{align*}
\|\phi(v)\|_{L^\infty([0,T];H^1)}&\leq  \|u_0\|_{H^1}+\left\|\int_0^t\e^{i(t-s)\Delta}[h(s)|x|^{-b}|v(s)|^{\alpha}v(s)]ds\right\|_{L^\infty([0,T];H^1)}\\
&\lesssim \|u_0\|_{H^1}+\|h|x|^{-b}|v|^{\alpha}v\|_{L^{q'}([0,T];L^{r',2})}+
\|h|x|^{-b}|v|^{\alpha}\nabla v\|_{L^{q'}([0,T];L^{r',2})}\\
&\vspace{1cm} +\|h|x|^{-b-1}|v|^{\alpha} v\|_{L^{q'}([0,T];L^{r',2})}=:\|u_0\|_{H^1}+L_1+L_2+L_3,
\end{align*}
where $(q,r)=(\frac{2}{b+(\alpha+1)\tau},\frac{2}{1-b-(\alpha+1) \tau})$ is an admissible pair. Indeed,  since $ 0<b<1 $ and by taking $ 0<\tau<\min\{\frac{1-b}{\alpha+1}, \frac{\alpha}{\alpha+1}\} $  we see that $ 1-b-(\alpha+1) \tau>0 $, then
$$0<\frac{1}{r}< \frac{1}{2}. $$ 
By H\"older's inequality in space, then in time, we obtain
\begin{align*}
L_1\lesssim \||x|^{-b}\|_{L^{\frac{2}{b},\infty}}\|h|v|^{\alpha}v\|_{L^{q'}([0,T]:L^{\eta,2})}\lesssim \|v\|^\alpha _{L^\infty([0,T]:L^{\tilde{   n_1},2})}\|hv\|_{L^{q'}([0,T]:L^{2,2})},
\end{align*}
where
$ \frac{1}{r'}=\frac{1}{\eta}+\frac{b}{2} $,\  
$ \frac{1}{\eta}=\frac{\alpha}{\tilde{n_1}}+\frac{1}{2}$ and
$ \frac{1}{\tilde{ n_1}}=\frac{(\alpha+1)\tau}{2\alpha}$, then we get 
$\tilde{ n_1}\in (2,\infty) $.
So we can use the Sobolev embedding $H^1  \hookrightarrow L^{\tilde{ n}_1,2}$, we get
$$ L_1\lesssim \|h\|_{L^{q'}([0,T])}\|v\|^\alpha_{L^\infty([0,T]:H^1)} \|v\|_{L^\infty([0,T]:L^2)}.$$
By using the same lines as above,
we get
$$L_2\lesssim \|h\|_{L^{q'}([0,T])}\|v\|^\alpha_{L^\infty([0,T]:H^1)} \|\nabla v\|_{L^\infty([0,T]:L^2)}. $$
We continue to estimate $L_3$. By H\"older’s inequality in  Lorentz spaces, we obtain
\begin{align*}
L_3&\leq \||x|^{-b-1}\|_{L^{\frac{2}{b+1},\infty}}\|h|v|^\alpha v\|_{L^{q'}([0,T]:L^{\eta_2})}\\
&\lesssim \|v\|^\alpha_{L^{\infty}([0,T]:L^{\tilde{n_2},2})} \|hv\|_{L^{q'}([0,T]:L^{ \tilde{ n_2},2})}\\
&\lesssim \|h\|_{L^{q'}([0,T])}\|v\|^\alpha_{L^{\infty}([0,T]:L^{\tilde{ n_2},2})}\|v\|_{L^{\infty}([0,T]:L^{\tilde{ n_2},2})},
\end{align*}
where $ \frac{1}{r'}=\frac{1}{\eta_2}+\frac{b+1}{2} $,\  
$ \frac{1}{\eta_2}=\frac{\alpha}{\tilde{  n_2}}+\frac{1}{\tilde{n_2} }   $ and $\frac{1}{\tilde{  n_2}}=\frac{\tau}{2}, \ \tau \in (0,1)$.\\ 
We use the Sobolev embedding $H^1 \hookrightarrow  L^{\tilde{n}_2,2}$, we get
\begin{align*}
L_3\lesssim \|h\|_{L^{q'}([0,T])} \|v\|^{\alpha+1}_{L^\infty([0,T]:H^1)}.
\end{align*}
We thus get 
\begin{align*}
\|\phi(v)\|_{L^\infty([0,T]:H^1)}&\lesssim \|u_0\|_{H^1}+\|h\|_{L^{q'}([0,T])}\|v\|^{\alpha+1}_{L^\infty([0,T]:H^1)}\\
&\leq K\|u_0\|_{H^1}+K\|h\|_{L^{q'}([0,T])}M^{\alpha +1}.
\end{align*}
In addition, we have
\begin{align*}
\|\phi(v)-\phi(w)\|_{L^\infty([0,T]:L^2)}&\leq K\|h\|_{L^{q'}}M^{\alpha }\|v-w\|_{L^\infty([0,T]:L^2)}.
\end{align*}
Finally, We choose $ M>0 $ and $ T > 0  $ such that
$$
\left\{
    \begin{array}{l}
     K\|u_0\|_{H^1}+KM^{\alpha+1}\|h\|_{L^{q'}}\leq M, \\
      KM^{\alpha}\|h\|_{L^{q'}}<1 ,
    \end{array}
  \right.
 $$
  where  $ K $  is a constant depend on $ \alpha $, $ b $ and $ \mu $. \\
If $\alpha >0  $ we set $K\|u_0\|_{H^1}=\frac{M}{2} $
and we choose $T$ so that
\begin{equation}\label{NNN4}
K M^{\alpha}\|h\|_{L^{q'}([0,T])}\leq \frac{1}{2}.
\end{equation}
Note that
\begin{align*}
\|h\|^{q'}_{L^{q'}([0,T])}
     =\frac{1}{\alpha \Re e(a) q'}\left( 1-\e^{-\alpha \Re e(a)q' T}\right).
 \end{align*}
Therefore, if $ \Re e(a)< \frac{(2-b-(\alpha+1)\tau)(2^{\alpha+1}K^{\alpha+1}\|u_0\|_{H^1}^\alpha)^{\frac{2}{2-b-(\alpha+1)\tau}}}{2\alpha}  $,
then we 
choose $T>0  $ verifies:
$$ T_{\max}(u_0)
> \frac{2-b-(\alpha+1)\tau}{2\Re e(a) \alpha }\log\left(\frac{(2^{\alpha+1}K^{\alpha+1}\|u_0\|_{H^1}^\alpha)^{q'}}{((2^{\alpha+1}K^{\alpha+1}\|u_0\|_{H^1}^\alpha)^{q'}-\alpha \Re e(a) q'}\right). $$
 In addition,  if $ \Re e(a)\geq \frac{(2-b-(\alpha+1)\tau)(2^{\alpha+1}K^{\alpha+1}\|u_0\|_{H^1}^\alpha)^{\frac{2}{2-b-(\alpha+1)\tau}}}{2\alpha} $,
we have a global solution.
So the fixed point argument on $E_2$ with the above values of $M$ and $T$ works and gives  a unique solution $v\in E_2$ of (\ref{2}), as well we get $(ii)$.
 Moreover, for any admissible pair $(p,l) $, $v \in L^p([0,T]:W^{1,l}_2 )$,  so we get $(i) $.\\
  Since the solution is constructed on an interval depending on $\|u_0\|_{H^1}$ we deduce the blowup alternative.
\end{proof}
\section{Scattering }
This section is devoted to study the  effect of the damping term on the scattering.
We first give the following  lemma.
\begin{lemma}
\label{scatteringcreteria}
Let $ N\geq 3 $,  $\mu \in \C $,  $ 0<b<\min\{2,N-2\}$, $a\in \C$ such that $\Re e (a)>0$ and $0 <\alpha\leq \frac{4-2b}{N-2}$.  Let $u$ be a  global solution of \eqref{INLS} and $v$ defined by  \eqref{transformation}. 
Then the following hold.
\begin{itemize}
\item[(i)] If $v$ scatters in $H^1$ then $v\in L^q ((0,\infty);W^{1,r}_2)$ for any admissible pair $(q,r).$ 
\item[(ii)] Let 
$\max\left({1\over N},{N-2\over 2N}, {1\over \alpha+1}\left({1\over 2}-{b\over N}+{\alpha\over N}\right)\right)< {1\over r}\leq\min\left({1\over 2}, {1\over \alpha+1}\left({1\over 2}-{b\over N}+{\alpha+2\over N}\right)\right),$ and $q$ so that $(q,r)$ be an admissible pair or $(q,r)=(\infty,2)$. If $v\in L^q ((0,\infty);W^{1,r}_2)$ then $v$ scatters  in $H^1$.
\end{itemize}  
\end{lemma}
\begin{remark}
{\rm The previous result in Part (ii), shows that global Strichartz regularity for some pairs suffices to have scattering. This is not the case without the damping term, that is $a=0.$ See \cite{CW}.}
\end{remark}
\begin{proof} (i) We have by the Duhamel formula,
  $$v(t)=\e^{it\Delta}v_+-i\mu \int_t^\infty \e^{i(t-s)\Delta} \{h(s)|x|^{-b}|v|^{\alpha}v(s)\}ds.  $$
  We show first that $v\in L^\gamma((0,\infty); W^{1,\rho}_2).$ Using the Strichartz estimates, we obtain
  \begin{align*}
 \|v\|_{L^\gamma((\tau,t);W^{1,\rho}_2)}&\lesssim\|v_+\|_{H^1}+ \|\e^{-a\alpha s}\|_{L^{\theta}((\tau,t))}\|v\|_{L^\gamma((\tau,t);W^{1,\rho}_2)}^{\alpha+1}\\
   &\lesssim \|v_+\|_{H^1}+\e^{-a\alpha \tau}\|v\|_{L^\gamma((\tau,t);W^{1,\rho}_2)}^{\alpha+1}.
   \end{align*}
Choosing $\tau$ sufficiently large and by continuity argument, see \cite[Lemma 7.7.4, p. 244]{Cazenave}, we deduce that 
$$\|v\|_{L^\gamma((\tau,\infty);W^{1,\rho}_2)}\lesssim\|v_+\|_{H^1}.$$
The result on $(0,\tau)$ follows by local theory. The result for the other admissible pair follows by Strichartz estimates using the pair $(\gamma,\rho)$ on the right-hand side.

   (ii)  We first consider the case $\max\left({1\over N},{N-2\over 2N}, {1\over \alpha+1}\left({1\over 2}-{b\over N}+{\alpha\over N}\right)\right)< {1\over r}\leq\min\left({1\over 2}, {1\over \alpha+1}\left({1\over 2}-{b\over N}+{\alpha+2\over N}\right)\right).$ Let $${1\over r_1}= 1-\frac{b}{N}-{\alpha+1\over r}+\frac{\alpha}{N},\; \frac{1}{n_1}=\frac{1}{r}-\frac{1}{N}.$$
It is clear that $2\leq r_1\leq {2N\over N-2}$ and $ 2 \leq  n_1<\infty.$ Let $q_1$ be such that $(q_1,r_1)$ is an admissible pair. 
Using the Strichartz estimates, we get
   \begin{eqnarray*}
\|\e^{-it\Delta}v(t)-\e^{-it'\Delta}v(t')\|_{H^1}&\lesssim&  \|h |x|^{-b}|v|^\alpha v\|_{L^{q_1'} ((t,t');W^{1,r_1'}_2)}.
\end{eqnarray*}
We have $$\frac{1}{r_1'}= \frac{1}{r}+\frac{b}{N}+{\alpha\over n_1},\ 
  \frac{1}{q_1'}=\frac{\alpha}{q}+\frac{1}{q}+\frac{1}{\theta}, 
 $$ then by the  H\"older inequality in Lorentz spaces we have
$$\|h|x|^{-b}|v|^{\alpha}v\|_{L^{{q_1'}}((t,t');L^{r_1',2})}
 \lesssim \|h\|_{L^{\theta}([t,t'])} \|v\|^\alpha _{L^{{q}}(t,t'):L^{n_1,2})}\|v\|_{L^q((t,t');L^{r,2})}.
$$
We use the homogenous Sobolev embedding $\dot{W}^{1,r}_2  \hookrightarrow L^{n_1,2}$ to get
$$ \|h|x|^{-b}|v|^{\alpha} v\|_{L^{q_1'}((t,t');L^{r_1',2})}\lesssim \|h\|_{L^{\theta}([t,t'])}\|v\|^\alpha_{L^{q}((t,t'));\dot{W}_2^{1,r})}\|v\|_{L^q((t,t');L^{r,2})}.$$
Similarly, we have 
$$ \|h|x|^{-b}|v|^{\alpha} \nabla v\|_{L^{q_1'}((t,t');L^{r_1',2})}\lesssim \|h\|_{L^{\theta}([t,t'])}\|v\|^{\alpha+1}_{L^{q}((t,t'));\dot{W}_2^{1,r})}.$$
We also have 
$$ \|h|x|^{-b-1}|v|^{\alpha} v\|_{L^{q_1'}((t,t');L^{r_1',2})}\lesssim \|h\|_{L^{\theta}([t,t'])}\|v\|^{\alpha+1}_{L^{q}((t,t'));\dot{W}_2^{1,r})}.$$
We then get 
$$\|\e^{-it\Delta}v(t)-\e^{-it'\Delta}v(t')\|_{H^1}\lesssim  \|h\|_{L^{\theta}([t,t'])}\|v\|^{\alpha+1}_{L^{q}((t,t'));W_2^{1,r})}.$$
This shows that $v$ scatters.

We now consider the case $(q,r)=(\infty,2)$. It is clear that $(\gamma,\rho)$ satisfies hypotheses the first part (ii). So, to prove the result it suffices to show   that $\|v\|_{L^\gamma ((0,\infty);W^{1,\rho}_2)} $ is a finite quantity.  

Using the  Strichartz estimates, we have for any  $0<\tau <T<\infty$,
\begin{align*}
\|v\|_{L^\gamma((\tau,T);W^{1,\rho}_2)}&\leq \|\e^{it\Delta}v(\tau)\|_{L^\gamma((\tau,T);W^{1,\rho}_2)}+\left\|\int_{\tau}^t \e^{i(t-s)\Delta}[h(s)|x|^{-b}|v(s)|^\alpha v(s)]ds\right\|_{L^\gamma((\tau,T);W^{1,\rho}_2)}\\
&\lesssim \|v(\tau)\|_{H^1}+\|h|x|^{-b}|v|^\alpha v\|_{L^{\gamma'}((\tau,T);W^{1,\rho'}_2)}\\
&\lesssim  \|v\|_{L^\infty([0,\infty);H^1)}+\|h|x|^{-b}|v|^{\alpha} v\|_{L^{\gamma'}((\tau,T);L^{\rho',2})}+\|h|x|^{-b}|v|^{\alpha}\nabla v\|_{L^{\gamma'}((\tau,T);L^{\rho',2})}\\
& \vspace{1 cm}+\|h|x|^{-b-1}|v|^{\alpha
}v\|_{L^{\gamma'}((\tau,T);L^{\rho',2})}:=\|v\|_{L^\infty([0,\infty);H^1)}+I+{II}+{III}.
\end{align*}
By Lemma~\ref{estimationstri}, we get 
$${I}\lesssim \|h\|_{{L^{\theta}}((\tau,T))}\|\nabla v\|^\alpha_{L^{\gamma}((\tau,T);L^{\rho,2})} \|v\|_{L^{\gamma}((\tau,T);L^{\rho,2})},$$

$$ {II}\lesssim \|h\|_{L^{\theta}((\tau,T))}\|\nabla v\|^{\alpha +1}_{L^{\gamma}((\tau,T);L^{\rho,2})}$$
 
 and 
 
 \begin{align*}
{III}
&\lesssim \|h\|_{L^{\theta}((\tau,T))}\|\nabla v\|^{\alpha+1}_{L^{\gamma}((\tau,T);L^{\rho,2})}.
 \end{align*}

Thus we obtain
\begin{align*}
\|v\|_{L^\gamma((\tau,T);W^{1,\rho}_2)}&\lesssim \|v\|_{L^\infty([0,\infty);H^1)}+\|h\|_{L^{\theta}((\tau,T))}\|v\|^{\alpha+1}_{L^{\gamma}((\tau,T);W^{1,\rho}_2)}.
\end{align*}

For large $\tau >0$, we get  $\|v\|_{L^\gamma((\tau,T);W^{1,\rho}_2)}\lesssim \|v\|_{L^\infty([0,\infty);H^1)}$. Fix such $\tau$, since  $T$ is arbitrary, we obtain $\|v\|_{L^\gamma((\tau,\infty);W^{1,\rho}_2)}\lesssim \|v\|_{L^\infty([0,\infty);H^1)}$. This implies that $\|v\|_{L^\gamma((\tau,\infty);W^{1,\rho}_2)}<\infty.$
Since  $v\in  L^\gamma_{loc}([0,\infty);W^{1,\rho}_2) $,
 then   $v\in L^\gamma((0,\infty);W^{1,\rho}_2) $. This completes the proof of the lemma.
\end{proof}

\begin{proof}[Proof of Theorem \ref{scatteringforlargea}]
We give the proof for $s=1.$ The case $s=0$ follows similarly. By Lemma \ref{scatteringcreteria}, it suffices to show that $v$ is global and that  $v\in L^\gamma((0,\infty); W^{1,\rho}_2).$ 

(i) It follows by the proof of the local existence that we have $T_{\max}>T_0(\|v_0\|_{H^s})>0$ with $T_0$ independent of $a.$ Using the  Strichartz estimates, we have for any  $T \in (T_0, T_{\max})$,
\begin{align*}
\|v\|_{L^\gamma((T_0,T):W^{s,\rho}_2)}&\leq \|\e^{i(t-T_0)\Delta}v(T_0)\|_{L^\gamma((T_0,T):W^{s,\rho}_2)}+\\ & \left\|\int_{T_0}^t \e^{i(t-\sigma)\Delta}[h(\sigma)|x|^{-b}|v(\sigma)|^{\alpha}v(\sigma)]d\sigma\right\|_{L^\gamma((T_0,T):W^{s,\rho}_2)}\\
&\lesssim \|v(T_0)\|_{H^s}+\|h|x|^{-b}|v|^{\alpha}v\|_{L^{\gamma'}((T_0,T):W^{s,\rho'}_2)}\\
&\lesssim  \|v(T_0)\|_{H^s}+\|h|x|^{-b}|v|^{\alpha}v\|_{L^{\gamma'}((T_0,T):L^{\rho',2})}+\|h|x|^{-b-s}|v|^{\alpha}v\|_{L^{\gamma'}((T_0,T):L^{\rho',2})}\\
&\vspace{1 cm}+\|h|x|^{-b}|v|^{\alpha}(-\Delta)^{\frac{s}{2}} v\|_{L^{\gamma'}((T_0,T):L^{\rho',2})}.
\end{align*}
Next, we use the inequalities  (\ref{j1}) and (\ref{j2}) given by Lemma~\ref{estimationstri}, we get
\begin{align*}
\|h|x|^{-b}|v|^{\alpha}v\|_{L^{\gamma'}((T_0,T):L^{\rho',2})}\lesssim { \e^{-\alpha\theta \Re e(a)T_0}\over \alpha\theta \Re e(a)} \|(-\Delta)^{\frac{s}{2}}v\|^{\alpha}_{L^\gamma((T_0,T):L^{\rho,2})}
\|v\|_{L^\gamma((T_0,T):L^{\rho,2})},
\end{align*}
$$ \|h|x|^{-b}|v|^{\alpha}(-\Delta)^{\frac{s}{2}} v\|_{L^{\gamma'}((T_0,T):L^{\rho',2})}\lesssim  { \e^{-\alpha\theta \Re e(a)T_0}\over \alpha\theta \Re e(a)}\|(-\Delta)^{\frac{s}{2}}v\|^{\alpha+1}_{L^\gamma((T_0,T):L^{\rho,2})}, $$
and
\begin{align*}
\|h|x|^{-b-s}|v|^{\alpha}v\|_{L^{\gamma'}((T_0,T):L^{\rho',2})}&\lesssim { \e^{-\alpha\theta \Re e(a)T_0}\over \alpha\theta \Re e(a)}\|(-\Delta)^{\frac{s}{2}}v\|^{\alpha +1}_{L^\gamma((T_0,T):L^{\rho,2})}.
\end{align*}
Using the local construction of the solution and the Strichartz estimates, there exists $C>0$ independent of $a$ such that  $\|v(T_0)\|_{H^s}\leq C\|v_0\|_{H^s}.$  Then we get,
$$\|v\|_{L^\gamma((T_0,T):W^{s,\rho}_2)}\leq C\|v_0\|_{H^s}+C{ \e^{-\alpha\theta \Re e(a)T_0}\over \alpha\theta \Re e(a)}\|v\|^{\alpha +1}_{L^\gamma((T_0,T):W^{s,\rho}_2)}.$$

By taking $ \Re e(a)\geq a^*>0 $, where $a^*$ satisfies $${ \e^{-\alpha\theta a^*T_0}\over \alpha\theta a^*}={1\over (2C)^{\alpha+1}\|v_0\|_{H^s}^\alpha},$$ the continuity argument implies that for any $T \in (T_0, T_{\max})$,
$$\|v\|_{L^\gamma((T_0,T):W^{s,\rho}_2)}\leq \tilde{C}\|v_0\|_{H^s}.$$
where the constant $ \tilde{C}>0 $ is independent of $ T $. Then we obtain
$$\|v\|_{L^\gamma((T_0,T_{\max}):W^{s,\rho}_2)}\lesssim  \tilde{C}\|v_0\|_{H^s}.$$
Since $ v \in L^\gamma_{loc}([0,T_{\max}):W^{s,\rho}_2)  $, then  for all $ \Re e(a)>a^* $,
\begin{equation}\label{equ}
 \|v\|_{L^\gamma((0,T_{\max}):W^{s,\rho}_2)}\leq C(\|v_0\|_{H^s}).
 \end{equation}
 which, by the blow-up alternative, implies $T_{\max}(u_0)=\infty $ and $v$ scatters.

(ii) Let $T_0(u_0)$ be given by Theorem \ref{ge} (v). Using similar calculations as above, we have for any  $T \in (T_0, T_{\max})$,

$$ \|v\|_{L^\gamma((T_0,T):W^{s,\rho}_2)}\lesssim \|v(T_0)\|_{H^s}+ \e^{-\alpha \Re e(a)T_0}\|v\|^{\alpha+1}_{L^\gamma((T_0,T):{W}^{s,\rho}_2)}. $$
Using the local construction of the solution and the Strichartz estimates, there exists $C>0$ independent of $a$ such that  $\|v(T_0)\|_{H^s}\leq C.$ 
Then we get,
$$\|v\|_{L^\gamma((T_0,T):W^{s,\rho}_2)}\leq C+C\e^{-\alpha \Re e(a)T_0}\|v\|^{\alpha +1}_{L^\gamma((T_0,T):W^{s,\rho}_2)}.$$

By taking $ \Re e(a)\geq a^*(u_0)>0 $, where $a^*$ satisfies $$ \e^{-\alpha a^*T_0}={1\over (2C)^{\alpha+1}},$$ the continuity argument implies that for any $T \in (T_0, T_{\max})$,
$$\|v\|_{L^\gamma((T_0,T):W^{s,\rho}_2)}\leq \tilde{C}.$$
where the constant $ \tilde{C}>0 $ is independent of $ T $. Then we obtain
$$\|v\|_{L^\gamma((T_0,T_{\max}):W^{s,\rho}_2)}\lesssim  \tilde{C}.$$
The rest of the proof follows similarly as for the above case. This completes the proof of Theorem \ref{scatteringforlargea}.
\end{proof}

To prove Theorem ~\ref{sca1}, we need the following  functional (see  \cite{I}) 
\begin{align}
& H(t)=e^{2\Re\textit{e}(a)t} E(u(t))-\mu \frac{\Re\textit{e}(a)\alpha}{\alpha +2}\int_0^t  \e^{2\Re\textit{e}(a)s} \||x|^{-b} |u(s)|^{\alpha +2}\|_{L^1} ds.\label{hu}
\end{align}
By (\ref{energy}) we have  
$\partial_tH(t)=0$, that is $H$ is  conserved.
\begin{proof}[Proof of Theorem~\ref{sca1} ]
\textbf{First case}. Let $\mu <0 $ and $0 <\alpha< \frac{4-2b}{N-2}$. Since $H(t)$ is a  conserved quantity and $\mu<0$, we get
$$ \e^{2\Re\textit{e}(a)t} \|\nabla u(t)\|_{L^2}^2 \leq 2H(t)=2H(0)<\infty .$$
By the blow up alternative and Lemma~\ref{scatteringcreteria} (ii), we deduce that the solution is global and scatters. The first identity in Part 1) is proved above.
Applying now  the definition of $H(t)$  given by  (\ref{hu}), we see that
$$ 
E(u(t))\leq \e^{-2\Re e(a)t }H(0)$$
and
\begin{equation}\label{limite}
    \int_0^t \e^{2\Re\textit{e}(a)s} \||x|^{-b} |u(s)|^{\alpha +2}\|_{L^1} ds\leq\frac{\alpha +2}{|\mu|\Re e(a)\alpha} E(u_0).
\end{equation}
We may let $t\rightarrow \infty$ in (\ref{limite}) and  we obtain
$$\int_0^\infty \e^{2\Re\textit{e}(a)s} \||x|^{-b} |u(s)|^{\alpha +2}\|_{L^1} ds\leq \frac{\alpha +2}{|\mu|\Re e(a)\alpha} E(u_0). $$ This completes the proof of Part 1).

\textbf{Second case}: Let $ \mu >0 $ and  $0 <\alpha< \frac{4-2b}{N}$. 
By H\"older's inequality, we have
\begin{align*}
\frac{1}{2}\|\nabla v(t)\|^2_{L^2}&=H(t)+ \mu\frac{\e^{-\alpha \Re \textit{e}(a)t}}{\alpha +2}\||x|^{-b}|v(t)|^{\alpha +2}\|_{L^1}+\mu\frac{\Re \textit{e}(a)\alpha}{\alpha +2}\int_0^t  \e^{-\alpha \Re \textit{e}(a)s} \||x|^{-b} |v(s)|^{\alpha +2}\|_{L^1} ds\\
& \lesssim H(0) +\mu\frac{\e^{-\alpha \Re \textit{e}(a)t}}{\alpha +2}\||x|^{-b}\|_{L^{\frac{N}{b},\infty}}\||v(t)|^{\alpha +2}\|_{L^{\frac{N}{N-b},1}
}+\mu\frac{\Re \textit{e}(a)\alpha}{\alpha +2}\times\\
&\hspace{1cm}\int_0^t  \e^{-\alpha \Re \textit{e}(a)s} \||x|^{-b}\|_{L^{\frac{N}{b},\infty}}\||v(s)|^{\alpha +2}\|_{L^{\frac{N}{N-b},1}} ds\\
&\lesssim H(0) +\mu\frac{\e^{-\alpha \Re \textit{e}(a)t}}{\alpha +2}\|v(t)\|^{\alpha +2}_{L^{\frac{N(\alpha +2)}{N-b},\alpha +2}} +\mu\frac{\Re \textit{e}(a)\alpha}{\alpha +2} \int_0^t \e^{-\alpha \Re \textit{e}(a)s} \|v(s)\|^{\alpha +2}_{L^{\frac{N(\alpha +2)}{N-b},\alpha +2}} ds.
\end{align*}
 By the Gagliardo-Nirenberg inequality and  Calder\'on inequality, it follows that \begin{align*}
\frac{1}{2}\|\nabla v(t)\|^2_{L^2}&\lesssim   H(0) +\mu\frac{\e^{-\alpha \Re \textit{e}(a)t}}{\alpha +2}\|\nabla v(t)\|_{L^{2,\alpha +2}}^{\frac{N\alpha +2b}{2}} \| v(t)\|_{L^{2,\alpha +2}}^{\frac{4-2b-(N-2)\alpha}{2}}\\
& +\mu\frac{\Re \textit{e}(a)\alpha}{\alpha +2}\int_0^t  \e^{-\alpha \Re \textit{e}(a)s} \|\nabla v(s)\|_{L^{2,\alpha +2}}^{\frac{N\alpha +2b}{2}} \| v(s)\|_{L^{2,\alpha +2}}^{\frac{4-2b-(N-2)\alpha}{2}}ds\\
& \lesssim  H(0) +\mu\frac{\e^{-\alpha \Re \textit{e}(a)t}}{\alpha +2}\|\nabla v(t)\|_{L^{2,2}}^{\frac{N\alpha +2b}{2}} \| v(t)\|_{L^{2,2}}^{\frac{4-2b-(N-2)\alpha}{2}}\\
&+\mu\frac{\Re \textit{e}(a)\alpha}{\alpha +2}\int_0^t  \e^{-\alpha \Re \textit{e}(a)s} \|\nabla v(s)\|_{L^{2,2}}^{\frac{N\alpha +2b}{2}} \| v(s)\|_{L^{2,2}}^{\frac{4-2b-(N-2)\alpha}{2}}ds.
 \end{align*}
Put $f(t)=\|\nabla v(t)\|_{L^2} $.  We get
\begin{eqnarray}\label{in4}
\frac{1}{2} f(t)^2 &\leq &  H(0)+C_1\frac{\e^{-\alpha \Re \textit{e}(a)t}}{\alpha +2} f(t)^{\frac{N\alpha +2b}{2}} \| u_0\|_{L^{2}}^{\frac{4-2b-(N-2)\alpha}{2}} \nonumber\\&& +C_1\frac{\Re \textit{e}(a)\alpha}{\alpha +2}\int_0^t  \e^{-\alpha \Re \textit{e}(a)s}f(s)^{\frac{N\alpha +2b}{2}} \|u_0\|_{L^{2}}^{\frac{4-2b-(N-2)\alpha}{2}}ds,
 \end{eqnarray}
 where $C_1=C_1(b,\alpha,N)>0$.  Using the inequality,
 \begin{align*}
  \frac{C_1}{\alpha +2}\e^{-\alpha \Re \textit{e}(a)t} f(t)^{\frac{N\alpha +2b}{2}} \| u_0\|_{L^{2}}^{\frac{4-2b-(N-2)\alpha}{2}} 
  &=\frac{C_1}{(\alpha +2)}f(t)^{\frac{N\alpha +2b}{2}}\e^{-\frac{\alpha}{2}\Re \textit{e}(a)t} \e^{-\frac{\alpha}{2}\Re \textit{e}(a)t}  \| u_0\|_{L^{2}}^{\frac{4-2b-(N-2)\alpha}{2}}\\
  &\leq \frac{1}{4}\e^{-\frac{2\alpha}{N\alpha+2b}\Re \textit{e}(a)t} f(t)^2 +CA \e^{-\frac{2\alpha}{4-2b-N\alpha }\Re \textit{e}(a)t}\\
  &\leq \frac{1}{4}\ f(t)^2 +CA,
   \end{align*}
   where $A= \| u_0\|_{L^{2}}^{\frac{2(4-2b-(N-2)\alpha)}{4-2b-N\alpha}}$, we obtain
\begin{align*}
\frac{1}{2} f(t)^2\leq& H(0)+\frac{1}{4}f(t)^2+CA \\
 &+\Re e(a)\alpha\left[\int_0^t \frac{1}{4} e^{-\frac{2\alpha}{N\alpha+2b}\Re \textit{e}(a)s} f(s)^{2} ds+  CA \int_0^t\e^{-\frac{2\alpha}{4-2b-N\alpha }\Re \textit{e}(a)s}ds  \right]\\
\leq& H(0)+\frac{1}{4}f(t)^2+CA +\Re e(a)\alpha\left(\int_0^t \frac{1}{4}  e^{-\frac{2\alpha}{N\alpha+2b}\Re \textit{e}(a)s} f(s)^{2} ds\right),
\end{align*}
and so
$$ \frac{1}{4} f(t)^2\leq  H(0)_++C A+\Re e(a)\alpha\left(\int_0^t \frac{1}{4}  e^{-\frac{2\alpha}{N\alpha+2b}\Re \textit{e}(a)s} f(s)^{2} ds\right).$$
Finally, it follows from the Gronwall inequality that
$$f(t)^2\leq 4H(0)_++C(\Re e(a),b,\alpha,N)A.$$
We now conclude by the blow-up alternative
 that the solution is global. Furthermore, by
Lemma~\ref{scatteringcreteria} (ii), the solution  scatters. 

It remains to show (\ref{b}). The last estimate gives 
\begin{equation}\label{grad}
\|\nabla v(t)\|_{L^2}^2\lesssim \|\nabla u_0\|_{L^2}^2+A. 
\end{equation}
By the transformation (\ref{transformation}), it follows easily that 
\begin{equation}\label{grau}
    \|\nabla u(t)\|_{L^2}^2\lesssim \e^{-2\Re\textit{e}(a)t}\left( \|\nabla u_0\|_{L^2}^2+\| u_0\|_{L^{2}}^{\frac{2(4-2b-(N-2)\alpha)}{4-2b-N\alpha}}\right). 
\end{equation}
This proves \eqref{b}.

 Next, by the conservation of $H$ we have 
\begin{equation}
\label{conservationH2}
 \frac{\mu\e^{2\Re\textit{e}(a)t}}{\alpha +2}\||x|^{-b}|u(t)|^{\alpha +2}\|_1+\frac{\mu\Re\textit{e}(a)\alpha}{\alpha +2}\int_0^t  \e^{2\Re\textit{e}(a)s} \||x|^{-b} |u(s)|^{\alpha +2}\|_1ds=\frac{1}{2}\e^{2\Re\textit{e}(a)t}\|\nabla u(t)\|^2_2-H(0),
\end{equation}
By (\ref{grau}), we get \eqref{c} and \eqref{d}.

 \textbf{Third case.} 
Let $\alpha = \frac{4-2b}{N}.$ We recall from \cite[Corollary 2.3, p. 287]{Genoud} the following Gagliardo-Nirenberg inequality
\begin{equation}
\label{GNN}
\||\cdot |^{-b}|g|^{\alpha +2}\|_{L^1}\leq {\alpha+2\over 2}\left({\| g\|_2\over \| Q\|_2}\right)^\alpha \|\nabla g\|_{L^2}^2 ,\; g\in H^1.
\end{equation}
We have, for $0\leq t<T_{\max}(u_0),$
\begin{eqnarray*}
\frac{1}{2}f(t)^2&=&H(0)+ \mu\frac{\e^{-\alpha \Re \textit{e}(a)t}}{\alpha +2}\||x|^{-b}|v(t)|^{\alpha +2}\|_{L^1}+\mu\frac{\Re \textit{e}(a)\alpha}{\alpha +2}\int_0^t  \e^{-\alpha \Re \textit{e}(a)s} \||x|^{-b} |v(s)|^{\alpha +2}\|_{L^1} ds\\  &\leq& H(0)+\frac{\mu}{2} f(t)^{2} \left({\| u_0\|_2\over \| Q\|_2}\right)^\alpha +\frac{\mu\Re\textit{e}(a)\alpha}{2} \left({\| u_0\|_2\over \| Q\|_2}\right)^\alpha \int_0^t  \e^{-\alpha \Re\textit{e}(a)s}f(s)^{2} ds.
\end{eqnarray*}
Therefore, under the hypothesis on $u_0$, $E(u_0)>0$ and we get
 $$\frac{1}{2}\left(1-\mu \left({\| u_0\|_2\over \| Q\|_2}\right)^\alpha\right)f(t)^2\leq   E(0)+\frac{\mu\Re\textit{e}(a)\alpha}{2} \left({\| u_0\|_2\over \| Q\|_2}\right)^\alpha \int_0^t  \e^{-\alpha \Re\textit{e}(a)s}f(s)^{2} ds.$$
 Finally, it follows from the Gronwall inequality, we deduce that $T_{\max}(u_0)=\infty$ and that 
 \begin{equation}\label{f}
  f(t)^2 \leq {2E(0)\| Q\|_2^\alpha \over \| Q\|_2^\alpha-\mu \| u_0\|_2^\alpha}\exp \left({\mu\|u_0\|_2^{\alpha}\over \| Q\|_2^\alpha-\mu \| u_0\|_2^\alpha}  \right).
   \end{equation}
  By Lemma~\ref{scatteringcreteria} (ii), the solution scatters. Using the transformation (\ref{transformation}), we get 
$$\e^{2\Re \textit{e}(a)t}\|\nabla u(t)\|_{L^2}^2\leq  {2E(0)\| Q\|_2^\alpha \over \| Q\|_2^\alpha-\mu \| u_0\|_2^\alpha}\exp \left({\mu\|u_0\|_2^{\alpha}\over \| Q\|_2^\alpha-\mu \| u_0\|_2^\alpha}  \right).$$
This proves \eqref{intg1}. Estimate \eqref{intg2} follows by  \eqref{conservationH2}.
\end{proof}

We  give now the proof of time decay
when the solution exponentially scatters.

\begin{proof}[Proof of Theorem \ref{decayratescattering}]
 By Duhamel's formula, we have
  $$v(t)=\e^{it\Delta}v_+-i\mu \int_t^\infty \e^{i(t-s)\Delta} \{h(s)|x|^{-b}|v|^{\alpha}v(s)\}ds.  $$
  Using the Strichartz estimates, we obtain
  \begin{align*}
  \|v(t)-\e^{it\Delta}v_+\|_{H^1}&=|\mu|\left\|\int_t^\infty \e^{i(t-s)\Delta} \{h(s)|x|^{-b}|v|^{\alpha}v\}ds\right\|_{H^1}\\
   &\lesssim \|\e^{-a\alpha s}\|_{L^{\theta}((t,\infty))}\|v\|_{L^\gamma((t,\infty):W^{1,\rho}_2)}^{\alpha+1}\\
   &\lesssim \e^{-a\alpha t}\|v\|_{L^\gamma((t,\infty):W^{1,\rho}_2)}^{\alpha+1}.
   \end{align*}
By the hypothesis $v$ scatters  then using Lemma \ref{scatteringcreteria} (i), $\|v\|_{L^\gamma((0,\infty):W^{1,\rho}_2)}<\infty$ hence $\|v\|_{L^\gamma((t,\infty);W^{1,\rho}_2)}
   \rightarrow 0\  \mbox{as}\  t\rightarrow \infty.$ This completes the proof of the theorem.
   \end{proof}
    We close this section by the following complementary result.
  \begin{proposition}\label{scattering}
Let $ N\geq 3 $,  $\mu >0$,  $ 0<b<\min\{2,N-2\}$, $a\in \C$ such that $\Re e (a)>0$ and $ \frac{4-2b}{N}<\alpha\leq \frac{4-2b}{N-2}$.  Let $u$ be a $H^1$ the global solution of \eqref{INLS}. Assume that one of the following holds.
\begin{itemize}
  \item[(i)] $ \displaystyle\liminf_{t\rightarrow \infty}\e^{-\Re\textit{e}(a)t}\|\nabla u(t)\|_{L^2(\R^N)}=0 $.
\item[(ii)] $ \displaystyle\lim_{t\rightarrow \infty}\e^{-\kappa \Re \textit{e}(a) t}\|\nabla u(t)\|_{L^2(\R^N)}=0 $, where $ \kappa =\frac{4-2b-(N-2)\alpha}{N\alpha-4+2b} $.
\end{itemize}  
Then $ u$ exponentially scatters.
\end{proposition}

To prove Proposition~\ref{scattering}, we need the following. 

\begin{lemma}\label{2.9}
Let $ N\geq 3 $,  $\mu >0$,  $ 0<b<\min\{2,N-2\}$, $a\in \C$ such that $\Re e (a)>0$ and $0 <\alpha< \frac{4-2b}{N-2}$.  Let $u$ be a $H^1$ global solution to
\eqref{intINLS}. Then, we have
$$\liminf_{t\rightarrow \infty}K(u(t)) \leq 0,$$
where $K$ is given by \eqref{energy1}.
\end{lemma}

\begin{proof}
The proof of the previous  lemma can be done similarly as in \cite[Lemma 2.9 page 770]{I}. So we omit the details.  
\end{proof}

\begin{proof}[Proof of Proposition~\ref{scattering} ] 
 
Part (i) can be obtained  as in \cite[Theorem 1.4 (b), page 765]{I}, so we omit it. We now consider the case $(ii)$. Using the Gagliardo-Nirenberg, we get 
 \begin{align*}
K(t)&=\|\nabla u(t)\|_{L^2}^2-\mu \||x|^{-b}|u(t)|^{\alpha +2}\|_{L^{1,1}}\\
&\geq \|\nabla u(t)\|_{L^2}^2- \mu C\|\nabla u(t)\|_{L^{2,2}}^{\frac{N\alpha +2b}{2}}\|u(t)\|_{L^{2,2}}^{\frac{4-2b-(N-2)\alpha}{2}}\\
&\geq \|\nabla u(t)\|_{L^2}^2\left[ 1-C\mu \e^{-\frac{4-2b-(N-2)\alpha}{2}\Re \textit{e}(a)t}\|u_0\|_{L^{2}}^{\frac{4-2b-(N-2)\alpha}{2}}\|\nabla u(t)\|_{L^{2}}^{\frac{N\alpha -4+2b}{2}}\right],
\end{align*}
where $C$ is  the constant given by the Gagliardo-Nirenberg inequality. 
Since $ \displaystyle\lim_{t\rightarrow\infty}e^{-\kappa \Re\textit{e}(a)t }\|\nabla u(t)\|_{L^2}=0 $ then there exists a large $ T> 0  $ such that
$$ K(t)\geq \frac{1}{2}\|\nabla u(t)\|_{L^2}^2, $$
for all $ t > T $. By Lemma~\ref{2.9} we get 
$$ \liminf_{t\rightarrow \infty}\|\nabla u(t)\|_{L^2}^2 =0.$$
The result follows then by $(i)$.
\end{proof}
  \section{Blow up} 
In  this section we prove Theorem~\ref{blow}. We may assume without loss of generality that $\mu=1$.  Before proceeding to the proof, we need to introduce some definitions.
 Let $v$ be the solution of (\ref{2}), we set 
  \begin{align}
&{\mathcal{M}}(t):=\|v(t)\|_{L^2}^2,\\
     &{\mathcal{E}}(t):=\frac{1}{2}\|\nabla v(t)\|_{L^2}^2-\frac{\e^{-\alpha at}}{\alpha +2}\||x|^{-b}|v(t)|^{\alpha +2}\|_{L^1}\label{e}=e^{2at}E(u(t)),\\
   &H(t):=\frac{1}{2}\|\nabla v(t) \|^2_{L^2}-  \frac{\e^{-\alpha at}}{\alpha +2}\||x|^{-b}|v(t)|^{\alpha +2}\|_{L^1}- \frac{a\alpha}{\alpha +2}\int_0^t  \e^{-\alpha as} \||x|^{-b} |v(s)|^{\alpha +2}\|_{L^1} ds,\label{h}\\
   &{\mathcal{I}}(t):=\|xv(t)\|_{L^2}^2,\\
   &{\mathcal{V}}(t):=\Im \textit{m}\int_{\R^n} x.\nabla v(t,x).\overline{v(t,x)}dx,\\
   &P(t):=\frac{1}{2}\|\nabla v(t)\|_{L^2}^2-\e^{-\alpha at}\frac{N\alpha +2b}{4(\alpha +2)}\||x|^{-b}|v(t)|^{\alpha +2}\|_{L^1} \label{p}.
\end{align}

 We can check that

\begin{align}
&{\mathcal{M}}'(t)=0,\label{mm}\\
\label{ee}
    &  {\mathcal{E}}'(t)=\frac{a\alpha}{\alpha +2}\e^{-\alpha at}\||x|^{-b}|v(t)|^{\alpha +2}\|_{L^1},\\
    \label{ii}
    & {\mathcal{I}}'(t)=4 {\mathcal{V}}(t),\\
  \label{vv}
    & {\mathcal{V}}'(t)=4 P(t).
    \end{align}
  We know give the proof of the theorem.
\begin{proof}[Proof of Theorem~\ref{blow} ] Let $a>0 $ and recall that
    \begin{equation*}
  \gamma =\frac{4\alpha}{N\alpha -4+2b}a.
     \end{equation*}
   Then, since $\alpha>\frac{4-2b}{N}$, we have $\gamma >0$. By (\ref{ee}),  (\ref{e}) and (\ref{p}),  
   we have 
   \begin{align}
   \frac{d}{dt}\left(\e^{-\gamma t}{\mathcal{E}}(t)\right) &=-\gamma \e^{-\gamma t}{\mathcal{E}}(t)+\e^{-\gamma t}{\mathcal{E}}'(t)\nonumber\\
   &=-\gamma \e^{-\gamma t}\left( \frac{1}{2}\|\nabla v(t)\|_{L^2}^2-\frac{\e^{-\alpha at}}{\alpha +2}\||x|^{-b}|v(t)|^{\alpha +2}\|_{L^1}-\frac{1}{\gamma} \e^{-\alpha at}\frac{a\alpha}{\alpha +2}\||x|^{-b}|v(t)|^{\alpha +2}\|_{L^1} \right)\nonumber\\
  &= -\gamma \e^{-\gamma t}\left( \frac{1}{2}\|\nabla v(t)\|_{L^2}^2- \e^{-\alpha at}\frac{1}{\alpha +2}(1+\frac{N\alpha -4+2b}{4})\||x|^{-b}|v(t)|^{\alpha +2}\|_{L^1}\right)\nonumber\\
  &=-\gamma \e^{-\gamma t}P(t).\label{eee}
   \end{align}
   Moreover, by (\ref{ii}) and (\ref{vv}), we have  
   \begin{equation}\label{iii}
   \frac{d}{dt}\left(\e^{-\gamma t}{\mathcal{I}}(t)\right)=\e^{-\gamma t}(-\gamma {\mathcal{I}}(t)+4{\mathcal{V}}(t)),
   \end{equation}
   \begin{equation}\label{vvv}
   \frac{d}{dt}\left(\e^{-\gamma t}{\mathcal{V}}(t)\right)=\e^{-\gamma t}(-\gamma {\mathcal{V}}(t)+4P(t)).
   \end{equation}
Using $P(t)\leq {\mathcal{E}}(t)$ and   (\ref{eee}), we get 
   \begin{equation}\label{in}
   \e^{-\gamma t}P(t) \leq \e^{-\gamma t}{\mathcal{E}}(t)=E_0-\gamma \int_0^t  \e^{-\gamma s} P(s)ds.
   \end{equation}
   Set
    $$\tilde{P}(t)=\int_0^t \e^{-\gamma s}P(s)ds, \ \tilde{V}(t)=\int_0^t \e^{-\gamma s} {\mathcal{V}}(s)ds. $$
   Then, we write (\ref{in}) as follows 
   $$\tilde{P}'(t) +\gamma\tilde{P}(t)\leq E_0.$$
   Since 
   $\tilde{P}(0)=0$, we have 
   \begin{equation}\label{nn}
   \e^{\gamma t}\tilde{P}(t)\leq \frac{\e^{\gamma t}-1}{\gamma}E_0.
   \end{equation}
   Integrating (\ref{vvv}), we get 
   $$\tilde{V}'(t)+\gamma \tilde{V}(t)=V_0+4\tilde{P}(t), $$
   and by (\ref{nn}), we have 
   $$\frac{d}{dt}\left(\e^{\gamma t}\tilde{V}(t)\right)=V_0\e^{\gamma t}+4 \e^{\gamma t} \tilde{P}(t)\leq V_0\e^{\gamma t}+\frac{4E_0}{\gamma} (\e^{\gamma t}-1).$$
   Since $ \tilde{V}(0)=0 $, we obtain
   \begin{equation}\label{in2}
   \e^{\gamma t} \tilde{V}(t)\leq \frac{V_0}{\gamma}(\e^{\gamma t}-1)+\frac{4E_0}{\gamma^2}(\e^{\gamma t}-1-\gamma t). 
   \end{equation}
   Moreover,  since ${\mathcal{I}}(t)\geq 0 $, by integrating (\ref{iii}) and using (\ref{in2}), we get 
   \begin{equation}\label{12}
      {\mathcal{I}}(t)=\e^{\gamma t}\e^{-\gamma t}{\mathcal{I}}(t)\leq  \e^{\gamma t}\{I_0+4\tilde{V}(t)\}\leq g(\gamma , t),
   \end{equation}
   where
   $$ g(\gamma , t)=I_0\e^{\gamma t}+\frac{4V_0}{\gamma}(\e^{\gamma t}-1)+\frac{16E_0}{\gamma^2}(\e^{\gamma t}-1-\gamma t). $$
   
 We write 
   $$ g(\gamma , t)={1\over \gamma^2} ({C_1}\e^{\gamma t}-C_2t-C_3), $$ 
    where we put
   $$C_1=I_0\gamma^2+4V_0\gamma+16E_0,\;  C_2=16E_0\gamma,\; C_3=4V_0\gamma+16E_0. $$
   
(i)  By the hypotheses on $\gamma$ we have $C_1\leq 0$ and 
   $$ g(\gamma , t)\leq {C_1\over 2}t^2+(I_0\gamma+4V_0)t+I_0. $$
    Since $t_0={2I_0\over \sqrt{(I_0\gamma+4V_0)^2-2I_0C_1}-(I_0\gamma+4V_0)}$ is the positive root of  ${C_1\over 2}t^2+(I_0\gamma+4V_0)t+I_0=0$. Thus, by (\ref{12}), we see that
 \begin{align}\label{R1}
  {\mathcal{I}}(t_0) \leq 0,
 \end{align}
and  the result follows. 

For the other cases, it suffices to find $t_0>0$ so that $ g(\gamma , t_0)\leq 0.$ We argue as follows:
   
   (ii) We have 
   $g(\gamma , t)=(I_0+\frac{4V_0}{\gamma})\e^{\gamma t}-\frac{4V_0}{\gamma}$ and $g(\gamma , t_0)=0$ for $t_0={1\over \gamma}\log\left({4V_0\over \gamma I_0+4V_0}\right).$
   
   (iii) By the hypotheses we have  $C_1,\; C_2>0$ and $t_0={1\over \gamma}\log({C_2\over \gamma C_1})>0$ verifies ${\partial g(\gamma,t_0)\over \partial t}=0$ and $g(\gamma , t_0)<0.$
   
   (iv) We have $C_1>0$ and $\gamma<\displaystyle {2\over  I_0}\left(|V_0|-{\sqrt{V^2_0-4E_0I_0}}\right)<\displaystyle {4\over I_0|V_0|}\left(V_0^2-2E_0I_0\right).$ Then the proof follows similarly to (iii).
   
   (v) By hypotheses we have $C_3<C_1<0<C_2.$ Then $g(\gamma , t)\leq {1\over \gamma^2} ({C_1}\e^{\gamma t}-C_3)$. This  completes the proof of the theorem.

\end{proof}

\end{document}